\documentclass [leqno, 11pt] {amsart}

\usepackage{amssymb, amsmath, amsfonts, amsthm, amsbsy, amscd}

\newtheorem{theorem}{Theorem}[section]
\newtheorem{proposition}{Proposition}[section]
\newtheorem{corollary}{Corollary}[section]
\newtheorem{lemma}{Lemma}[section]

\newtheorem{definition}{Definition}[section]

\theoremstyle{definition}
\newtheorem{example}{Example}[section]
\newtheorem{remark}{Remark}[section]

\numberwithin{equation}{section}

\newcommand{\Hpull}{\pi^{\ast}(H)}

\newcommand{\pplus}{\mathfrak{p}^{+}}

\newcommand{\xxbbbdb}[7]{
\begin{picture}(66,14)
\put(3,3){\line(1,0){60}}%
\put(3,3){\makebox(0,0){$\times$}}
\put(23,3){\makebox(0,0){$\times$}}%
\put(43,3){\makebox(0,0){$\bullet$}}%
\put(63,3){\makebox(0,0){$\bullet$}}%
\put(73,3){\makebox(0,0){\dots}}%
\put(83,3){\makebox(0,0){$\bullet$}}%
\put(83,2){\line(1,0){20}}%
\put(83,4){\line(1,0){20}}%
\put(93,3){\makebox(0,0){\scriptsize$\langle$}}%
\put(103,3){\makebox(0,0){$\bullet$}}%
\put(3,12){\makebox(0,0){\scriptsize $#1$}}%
\put(23,12){\makebox(0,0){\scriptsize $#2$}}%
\put(43,12){\makebox(0,0){\scriptsize $#3$}}%
\put(63,12){\makebox(0,0){\scriptsize $#4$}}%
\put(73,12){\makebox(0,0){\scriptsize $#5$}}%
\put(83,12){\makebox(0,0){\scriptsize $#6$}}%
\put(103,12){\makebox(0,0){\scriptsize $#7$}}\end{picture}
}

\newcommand{\xxbbbdbdual}[7]{
\begin{picture}(66,14)
\put(3,3){\line(1,0){60}}%
\put(83,2){\line(1,0){20}}%
\put(83,4){\line(1,0){20}}%
\put(93,3){\makebox(0,0){$\scriptsize\langle$}}%
\put(73,3){\makebox(0,0){\dots}}%
\put(-7,3){\makebox(0,0){$($}}
\put(3,3){\makebox(0,0){$\times$}}
\put(23,3){\makebox(0,0){$\times$}}%
\put(43,3){\makebox(0,0){$\bullet$}}%
\put(63,3){\makebox(0,0){$\bullet$}}%
\put(83,3){\makebox(0,0){$\bullet$}}%
\put(103,3){\makebox(0,0){$\bullet$}}%
\put(115,3){\makebox(0,0){$)^{\ast}$}}%
\put(140,3){\makebox(0,0){$\simeq$}}%
\put(3,12){\makebox(0,0){\scriptsize $#1$}}%
\put(23,12){\makebox(0,0){\scriptsize $#2$}}%
\put(43,12){\makebox(0,0){\scriptsize $#3$}}%
\put(63,12){\makebox(0,0){\scriptsize $#4$}}%
\put(73,12){\makebox(0,0){\scriptsize $#5$}}%
\put(83,12){\makebox(0,0){\scriptsize $#6$}}%
\put(103,12){\makebox(0,0){\scriptsize $#7$}}\end{picture}
}

\newcommand{\xbbbbdb}[7]{
\begin{picture}(66,14)
\put(3,3){\line(1,0){60}}%
\put(3,3){\makebox(0,0){$\times$}}
\put(23,3){\makebox(0,0){$\bullet$}}%
\put(43,3){\makebox(0,0){$\bullet$}}%
\put(63,3){\makebox(0,0){$\bullet$}}%
\put(73,3){\makebox(0,0){\dots}}%
\put(83,3){\makebox(0,0){$\bullet$}}%
\put(83,2){\line(1,0){20}}%
\put(83,4){\line(1,0){20}}%
\put(93,3){\makebox(0,0){\scriptsize$\langle$}}%
\put(103,3){\makebox(0,0){$\bullet$}}%
\put(3,12){\makebox(0,0){\scriptsize $#1$}}%
\put(23,12){\makebox(0,0){\scriptsize $#2$}}%
\put(43,12){\makebox(0,0){\scriptsize $#3$}}%
\put(63,12){\makebox(0,0){\scriptsize $#4$}}%
\put(73,12){\makebox(0,0){\scriptsize $#5$}}%
\put(83,12){\makebox(0,0){\scriptsize $#6$}}%
\put(103,12){\makebox(0,0){\scriptsize $#7$}}\end{picture}
}

\newcommand{\xbbbbdbdual}[7]{
\begin{picture}(66,14)
\put(3,3){\line(1,0){60}}%
\put(83,2){\line(1,0){20}}%
\put(83,4){\line(1,0){20}}%
\put(93,3){\makebox(0,0){\scriptsize$\langle$}}%
\put(73,3){\makebox(0,0){\dots}}%
\put(-7,3){\makebox(0,0){$($}}
\put(3,3){\makebox(0,0){$\times$}}
\put(23,3){\makebox(0,0){$\bullet$}}%
\put(43,3){\makebox(0,0){$\bullet$}}%
\put(63,3){\makebox(0,0){$\bullet$}}%
\put(83,3){\makebox(0,0){$\bullet$}}%
\put(103,3){\makebox(0,0){$\bullet$}}%
\put(115,3){\makebox(0,0){$)^{\ast}$}}%
\put(140,3){\makebox(0,0){$\simeq$}}%
\put(3,12){\makebox(0,0){\scriptsize $#1$}}%
\put(23,12){\makebox(0,0){\scriptsize $#2$}}%
\put(43,12){\makebox(0,0){\scriptsize $#3$}}%
\put(63,12){\makebox(0,0){\scriptsize $#4$}}%
\put(73,12){\makebox(0,0){\scriptsize $#5$}}%
\put(83,12){\makebox(0,0){\scriptsize $#6$}}%
\put(103,12){\makebox(0,0){\scriptsize $#7$}}\end{picture}
}
\newcommand{\bxbbbdb}[7]{
\begin{picture}(66,14)
\put(3,3){\line(1,0){60}}%
\put(3,3){\makebox(0,0){$\bullet$}}
\put(23,3){\makebox(0,0){$\times$}}%
\put(43,3){\makebox(0,0){$\bullet$}}%
\put(63,3){\makebox(0,0){$\bullet$}}%
\put(73,3){\makebox(0,0){\dots}}%
\put(83,3){\makebox(0,0){$\bullet$}}%
\put(83,2){\line(1,0){20}}%
\put(83,4){\line(1,0){20}}%
\put(93,3){\makebox(0,0){\scriptsize$\langle$}}%
\put(103,3){\makebox(0,0){$\bullet$}}%
\put(3,12){\makebox(0,0){\scriptsize $#1$}}%
\put(23,12){\makebox(0,0){\scriptsize $#2$}}%
\put(43,12){\makebox(0,0){\scriptsize $#3$}}%
\put(63,12){\makebox(0,0){\scriptsize $#4$}}%
\put(73,12){\makebox(0,0){\scriptsize $#5$}}%
\put(83,12){\makebox(0,0){\scriptsize $#6$}}%
\put(103,12){\makebox(0,0){\scriptsize $#7$}}\end{picture}
}

\newcommand{\bxbbbdbdual}[7]{
\begin{picture}(66,14)
\put(3,3){\line(1,0){60}}%
\put(83,2){\line(1,0){20}}%
\put(83,4){\line(1,0){20}}%
\put(93,3){\makebox(0,0){\scriptsize$\langle$}}%
\put(73,3){\makebox(0,0){\dots}}%
\put(-7,3){\makebox(0,0){$($}}
\put(3,3){\makebox(0,0){$\bullet$}}
\put(23,3){\makebox(0,0){$\times$}}%
\put(43,3){\makebox(0,0){$\bullet$}}%
\put(63,3){\makebox(0,0){$\bullet$}}%
\put(83,3){\makebox(0,0){$\bullet$}}%
\put(103,3){\makebox(0,0){$\bullet$}}%
\put(115,3){\makebox(0,0){$)^{\ast}$}}%
\put(140,3){\makebox(0,0){$\simeq$}}%
\put(3,12){\makebox(0,0){\scriptsize $#1$}}%
\put(23,12){\makebox(0,0){\scriptsize $#2$}}%
\put(43,12){\makebox(0,0){\scriptsize $#3$}}%
\put(63,12){\makebox(0,0){\scriptsize $#4$}}%
\put(73,12){\makebox(0,0){\scriptsize $#5$}}%
\put(83,12){\makebox(0,0){\scriptsize $#6$}}%
\put(103,12){\makebox(0,0){\scriptsize $#7$}}\end{picture}
}

\newcommand{\xxbbbdbt}[3]{
\begin{picture}(66,14)
\put(3,3){\line(1,0){20}}%
\put(23,2){\line(1,0){20}}%
\put(23,4){\line(1,0){20}}%
\put(33,3){\makebox(0,0){\scriptsize$\langle$}}%
\put(3,3){\makebox(0,0){$\times$}}
\put(23,3){\makebox(0,0){$\times$}}%
\put(43,3){\makebox(0,0){$\bullet$}}%
\put(3,12){\makebox(0,0){\scriptsize $#1$}}%
\put(23,12){\makebox(0,0){\scriptsize $#2$}}%
\put(43,12){\makebox(0,0){\scriptsize $#3$}}\end{picture}
}

\newcommand{\bxbbbdbt}[3]{
\begin{picture}(66,14)
\put(3,3){\line(1,0){20}}%
\put(23,2){\line(1,0){20}}%
\put(23,4){\line(1,0){20}}%
\put(33,3){\makebox(0,0){\scriptsize$\langle$}}%
\put(3,3){\makebox(0,0){$\bullet$}}
\put(23,3){\makebox(0,0){$\times$}}%
\put(43,3){\makebox(0,0){$\bullet$}}%
\put(3,12){\makebox(0,0){\scriptsize $#1$}}%
\put(23,12){\makebox(0,0){\scriptsize $#2$}}%
\put(43,12){\makebox(0,0){\scriptsize $#3$}}\end{picture}
}

\newcommand{\xx}[2]{
\begin{picture}(66,14)
\put(3,3){\makebox(0,0){$\times$}}%
\put(23,3){\makebox(0,0){$\times$}}%
\put(3,2){\line(1,0){20}}%
\put(3,4){\line(1,0){20}}%
\put(13,3){\makebox(0,0){\scriptsize$\langle$}}%
\put(3,12){\makebox(0,0){\scriptsize $#1$}}%
\put(23,12){\makebox(0,0){\scriptsize $#2$}}\end{picture}
}

\newcommand{\xxdual}[2]{
\begin{picture}(66,14)
\put(-7,3){\makebox(0,0){$($}}
\put(3,3){\makebox(0,0){$\times$}}%
\put(23,3){\makebox(0,0){$\times$}}%
\put(3,2){\line(1,0){20}}%
\put(3,4){\line(1,0){20}}%
\put(13,3){\makebox(0,0){\scriptsize$\langle$}}%
\put(3,12){\makebox(0,0){\scriptsize $#1$}}%
\put(23,12){\makebox(0,0){\scriptsize $#2$}}%
\put(35,3){\makebox(0,0){$)^{\ast}$}}%
\put(55,3){\makebox(0,0){$\simeq$}}%
\end{picture}
}

\newcommand{\kappah}{\kappa_{\mathcal{H}}}

\newcommand{\Det}{\text{Det}}
\newcommand{\sltwor}{\mathfrak{sl}(2, \mathbb{R})}
\newcommand{\cpl}{\mathcal{C}_{[L]}}
\newcommand{\chsub}{\text{Ch}}

\newcommand{\eno}{\text{End}}
\newcommand{\sgr}{SGr(2, \standrep)}

\newcommand{\cssp}{CSp(n, \mathbb{R})\cdot GL(2, \mathbb{R})}

\newcommand{\squat}{\mathbb{A}}
\newcommand{\lub}{\text{Lub}}

\newcommand{\dx}{\dot{x}}

\newcommand{\li}{\text{lie}}

\newcommand{\si}{\sigma}

\newcommand{\pr}{\partial}

\newcommand{\smp}{\mathfrak{sp}}
\newcommand{\spnr}{\mathfrak{sp}(n, \mathbb{R})}
\newcommand{\spl}{\mathfrak{sp}}
\newcommand{\grad}{\text{grad}}
\newcommand{\filt}{\text{filt}}
\newcommand{\gr}{\text{Gr}\,}

\newcommand{\integer}{\mathbb{Z}}
\newcommand{\en}{[\nabla]}

\newcommand{\ph}{\mathbb{P}(H)}

\newcommand{\paths}{\mathcal{S}}
\newcommand{\lie}{\mathfrak{L}}

\newcommand{\standu}{L}

\newcommand{\re}{\text{Re\,}}
\newcommand{\im}{\text{Im\,}}

\renewcommand{\xi}{\frac{\partial}{\partial x^{i}}}

\newcommand{\C}{\mathcal{C}}

\newcommand{\B}{\mathcal{B}}
\newcommand{\redframe}{\mathcal{G}_{0}}

\newcommand{\pframe}{\bar{\mathcal{G}}}
\newcommand{\ctractor}{\mathcal{P}}

\newcommand{\standrepz}{\mathbb{V}^{\times}}

\newcommand{\tmz}{T^{0}M}

\newcommand{\al}{\alpha}
\newcommand{\be}{\beta}
\newcommand{\ga}{\gamma}

\newcommand{\bg}{\boldsymbol{G}}
\newcommand{\btheta}{\boldsymbol{\Theta}}
\newcommand{\bomega}{\boldsymbol{\Omega}}
\newcommand{\bpsi}{\boldsymbol{\Psi}}
\newcommand{\levi}{\mathcal{L}}
\newcommand{\ztwo}{\mathbb{Z}^{2}}

\newcommand{\zk}{\mathbb{Z}^{k}}

\newcommand{\G}{\mathcal{G}}

\newcommand{\graded}{Gr\,}

\newcommand{\frameN}{\mathcal{F}}

\newcommand{\emf}{\mathcal{E}}

\newcommand{\adaptedframe}{\mathcal{G}}

\newcommand{\tractor}{\mathcal{T}}
\newcommand{\standrep}{\mathbb{V}}
\newcommand{\standw}{B}

\newcommand{\gl}{\mathfrak{gl}}
\newcommand{\sll}{\mathfrak{sl}}

\newcommand{\form}{\mathsf{K}}

\newcommand{\ct}{\Pi}
\newcommand{\rb}{T}
\newcommand{\proj}{\mathbb{P}}

\newcommand{\sym}{\text{Sym}\,}

\newcommand{\vect}{\text{Vec}}

\newcommand{\m}{\mathfrak{m}}
\newcommand{\q}{\mathfrak{q}}

\newcommand{\g}{\mathfrak{g}}
\newcommand{\f}{\mathfrak{f}}

\newcommand{\p}{\mathfrak{p}}

\newcommand{\ad}{ad}
\newcommand{\Ad}{Ad}

\newcommand{\tensor}{\otimes}

\newcommand{\Rea}{\mathbb R}
\newcommand{\rea}{\mathbb R}
\newcommand{\com}{\mathbb C}

\newcommand{\rpone}{\mathbb{P}^{1}(\mathbb{R})}

\newcommand{\tr}{\text{tr} \,}
\newcommand{\rank}{\text{rank\,}}

\newcommand{\ann}{Ann}
\newcommand{\K}{\mathcal{K}}

\newcommand{\ptm}{\mathbb{P}(TM)}
\newcommand{\piz}{\hat{\pi}}
\newcommand{\Hup}{\mathcal{H}}
\newcommand{\hup}{\mathcal{H}}

\begin{document}
\title{Contact Path Geometries}
\author{Daniel J. F. Fox} 
\address{School of Mathematics\\ Georgia Institute of Technology\\ 686 Cherry St.\\ Atlanta, GA 30332-0160, U.S.A.}
\email{fox@math.gatech.edu}
\date{\today}
\begin{abstract}
Contact path geometries are curved geometric structures on a contact manifold comprising smooth families of paths modeled on the family of all isotropic lines in the projectivization of a symplectic vector space. Locally such a structure is equivalent to the graphs in the space of independent and depedent variables of the family of solutions of a system of an odd number of second order ODE's subject to a single maximally non-integrable constraint. A subclass of contact path geometries is distinguished by the vanishing of an invariant contact torsion. For this subclass the equivalence problem is solved by constructing a normalized Cartan connection using the methods of Tanaka-Morimoto-\v{C}ap-Schichl. The geometric meaning of the contact torsion is described. If a secondary contact torsion vanishes then the locally defined space of contact paths admits a split quaternionic contact structure (analogous to the quaternionic contact structures studied by O. Biquard).
\end{abstract}
\thanks{During the preparation of this article the author benefited from visits to the Erwin Schr\"odinger Institute and the Universitat Aut\'onoma de Barcelona.}

\maketitle

\section{Introduction}
A \textbf{contact path geometry} is a smooth family of embedded one dimensional manifolds in a contact manifold, $(M, H)$, having the property that for each point, $L$, in an open subset of the projectivized contact hyperplane bundle, $\pi:\ph \to M$, there exists a unique path of the family passing through the basepoint $\pi(L) \in M$ and tangent to $L$. (See Definition \ref{nullpathdefn}). Two contact path geometries are regarded as equivalent if there is a contactomorphism mapping the paths of one family onto the paths of the other family. The flat model for the general contact path geometry is the family of all isotropic lines in the projectivization of a symplectic vector space (which has a canonical contact structure), and a contact path geometry is \textbf{locally flat} if it is locally equivalent to this model. The contact path geometries the contact paths of which are among the geodesics of some affine connection are called \textbf{contact projective structures} and were studied in \cite{Fox-cproj}.

Every contact path geometry on a three-dimensional contact manifold is locally equivalent to the graphs in the space with variables $(t, x, \dot{x})$ of the three-parameter family of solutions of a third order ordinary differential equation. The equivalence problem for third-order ordinary differential equations modulo contact transformations was studied first by K. W\"unschmann, \cite{Wunschmann}, and later solved by S.-S. Chern, \cite{Chern}. Here the discussion will focus on contact path geometries on contact manifolds of dimension $2n-1 \geq 5$. Locally such a contact path geometry is equivalent to a system of $2n-3$ second order ordinary differential equations depending on a parameter that is subject to a non-integrable constraint:
\begin{align}\label{contactpathodes}
&\ddot{x}^{i} = f^{i}(t, x^{j}, x^{0}, \dot{x}^{j}, \dot{x}^{0}, z), & &1 \leq i, j, p, q \leq 2n-4,\\
&\notag \ddot{x}^{0} = f^{0}(t, x^{j}, x^{0}, \dot{x}^{j}, \dot{x}^{0}, z) + \omega_{pq}\ddot{x}^{p}\dot{x}^{q},& &\omega_{ij} = - \omega_{ji},\\
&\notag  \dot{z} = x^{0} -t\dot{x}^{0} - \omega_{pq}x^{p}\dot{x}^{q},& &\Det(\omega_{ij}) \neq 0 & 
\end{align}
(Here $dz - x^{0}dt + tdx^{0} + \omega_{pq}x^{p}dx^{q}$ is a contact form). When a certain invariant, the \textbf{contact torsion}, vanishes, the first two of \eqref{contactpathodes} take the simpler forms
\begin{align}\label{cpathodenotorsion}
&\ddot{x}^{i} = -\tfrac{1}{3}\omega^{ip}\tfrac{\partial f}{\partial \dx^{p}} + \tfrac{1}{3}\dot{x}^{i}\tfrac{\partial f}{\partial \dx^{0}}, & &\ddot{x}^{0} = f -\tfrac{1}{3}\dot{x}^{p}\tfrac{\partial f}{\partial \dx^{p}},&
%& \dot{z} = x^{0} -t\dot{x}^{0} - \omega_{pq}x^{p}\dot{x}^{q}, & 
\end{align}
which evidently have a Hamiltonian flavor (see Remark \ref{torsionidentityremark} below).

From a Lie theoretic point of view, the flat model contact path geometry is an isotropic flag manifold of the form $G/P$, with $P \subset Sp(n, \rea) = G$ the parabolic subgroup corresponding to crossing the first two nodes on the Dynkin diagram $C_{n}$: 
\begin{center}
\parbox[r][1\totalheight][c]{158pt}{$\xxbbbdb{}{}{}{}{}{}{}$\\}
\end{center}

\noindent
When $n = 2$, these nodes are connected by a double bond and $P$ is a Borel subgroup, whereas when $n > 2$, they are connected by a single bond, and this accounts for the special nature of the three-dimensional case. The identity of the Dynkin diagrams $C_{2}$ and $B_{2}$ amounts to the special isomorphism $\smp(2, \rea) \simeq \mathfrak{so}(2, 3)$ induced by the action of the symplectic group on the space of trace-free two-forms on a $4$-dimensional symplectic vector space. This means that three-dimensional contact path geometries may be viewed as curved geometries modeled on a homogeneous quotient of the conformal Lie group, $SO(2, 3)$, and this was Chern's point of view in \cite{Chern}.

A solution of the equivalence problem for contact path geometries is taken here to mean the construction of a canonical Cartan connection on a principal bundle over the projectivized contact hyperplane bundle, $\ph$. Each contact path geometry of dimension at least five determines a $4$-step filtration of the tangent bundle, $T\ph$. A basic observation is that there exists a subclass of contact path geometries distinguished by the vanishing of the contact torsion, $\ct$, which is a tensorial object definable directly from the data determining a contact path geometry. The existence of this subclass is not obvious \textit{a priori}. The Lie bracket of vector fields induces a fiberwise algebraic Lie bracket on the associated graded of the filtered tangent bundle, called by N. Tanaka the symbola algebra of the filtration, and the symbola algebra of the filtration of $T\ph$ is isomorphic at every point of $\ph$ to that of the flat model if and only if $\ct = 0$. It is shown in Theorem \ref{cpathnormal} that in the event that $\ct = 0$ the general prolongation theory developed by N. Tanaka - T. Morimoto - A. \v{C}ap - H. Schichl can be applied to produce a canonical regular, normal Cartan connection. On the other hand, if the contact torsion does not vanish, the existing general prolongation procedures cannot be applied. Precisely, if for a contact path geometry there holds $\Pi \neq 0$ on an open subset $\B \subset \ph$, then a Cartan connection on $\B$ inducing the given contact path geometry can be neither regular nor normal; the first claim follows straightforwardly from the definition of regularity while the second claim, which is not obvious, is proved in Proposition \ref{normalimpliesctf} by analyzing harmonic curvature components. While it can be shown that there exist Cartan connections inducing a contact path geometry with non-vanishing contact torsion, there has not been identified a natural curvature normalization to impose on such a connection, and it remains a problem to solve the equivalence problem for contact path geometries in this generality. 

The space of contact paths in the flat model contact path geometry is parameterized by an isotropic Grassmannian of the form $Sp(n, \rea)/P$. The tangent bundle of this isotropic Grassmannian carries a filtration modeled on a split quaternionic Heisenberg algera, and it is interesting to consider when the (locally defined) space of contact paths of a contact torsion free contact path geometry admits such a \textbf{split quaternionic contact structure}. The split quaternionic contact structures that sometimes arise on the space of contact paths are interesting in their own right, being cousins of the quaternionic contact structures studied by O. Biquard in \cite{Biquard-big} in the context of asymptotically symmetric Einstein metrics. The two structures correspond respectively to homogeneous quotients of the split and compact real forms of $Sp(n, \com)$. A thorough discussion of split quaternionic contact structures would require more space. Here the prolongation machinery discussed above is used to show that isomorphism classes of split quaternionic contact structures correspond to isomorphism classes of regular, normal $(\mathfrak{sp}(n, \rea), P)$ Cartan connections, and the underlying data is identified in terms of tensors on a codimension three distribution. 

Contact path geometries are analogous in many respects to the structures studied by C. Lebrun in \cite{Lebrun}. Lebrun worked in the holomorphic category, but his results have the following smooth analogoues. Associate to a conformal class of pseudo-Riemannian metrics the projectivization of the null cone common to these metrics and consider \textbf{null path geometries}, namely families of paths in the base manifold such that for each one-dimensional subspace, $L$, in an open subset of the projectivized null cone there exists a unique path of the family tangent to $L$ at the basepoint of $L$. If these null paths are among the geodesics of some affine connection, the structure is a null projective structure. 
%THIS IS ONLY TRUE FOR NULL PATH GEOMETIES defined on entire null cone
%(in the holomorphic category this is automatic). 
Given a null projective structure there is associated to each representative metric a unique affine connection having among its geodesics the given null paths, making parallel the given metric, and satisfying certain natural normalizations on its torsion (these representative connections will be the Levi-Civita connections of the metric only in the event that the null projective structure under consideration is that null projective structure induced by the geodesics of the metrics). A piece of the torsion of the associated connection does not depend on the chosen metric, and this \textbf{conformal torsion} is formally analogous to the contact torsion of a contact projective structure. In particular, if it vanishes there exists on the space of null paths a canonical contact structure.

More generally, a conformal class of Lagrangians, $[L]$, is an equivalence class of smooth functions on the punctured tangent bundle any two of which differ by multiplication by a non-vanishing function on the base. When the image under the projection $\ptm \to M$ of the common zero locus, $\cpl$, in $\ptm$ of $[L]$, is $M$, it makes sense to speak of $\cpl$ null path geometries and projective structures. The conformal null path geometries are the special case in which representatives of $[L]$ are the fiberwise quadratic functions induced on the tangent bundle by the metrics representing a pseudo-Riemannian metric. Contact path geometries are the special case in which representatives of $[L]$ are the fiberwise linear functions induced on the tangent bundle by the contact one-forms. When each representative of $[L]$ is a non-degenerate Lagrangian, the generalized conformal structure $\cpl$ carries an even contact structure induced via $[L]$ from the contact structure on $\proj(T^{\ast}M)$, and so a characteristic line bundle. For example, a conformal null projective structure is torsion free if and only if its paths are the images in $M$ of the characteristics. From this point of view the contact path geometries are highly degenerate, in that the associated conformal class, $[L]$, comprises highly degenerate Lagrangians. In this case there is some structure on $\cpl = \ph$ analogous to that of the even contact structure on $\cpl$ for non-degenerate $[L]$. Precisely, there is a codimension $2$ distribution on $\ph$ which carries a family of symplectic structures, such that the skew complement of the line bundle generating the contact path geometry is independent of the choice of symplectic structure. The contact torsion vanishes if and only if this skew-complement is generated by Lie brackets of sections of the generating line bundle and certain vertical vector fields, and its characteristic system is the generating line bundle if and only if there vanishes an invariantly defined secondary contact torsion. %All this is explained in Section \ref{contacttorsion}. 

Using the theory of harmonic curvature components of a regular, normal Cartan connection, it is shown that when there vanishes the secondary contact torsion the resulting multicontact structure on the space of contact paths is in fact a split quaternionic contact structure. This is analogous to D. Grossman's result that a projective structure which is torsion free when viewed as a path geometry is flat, and both results follow by applying A. \v{C}ap's theory of correspondence spaces for parabolic geometries, \cite{Cap}. 

There is a correspondence realizing three-dimensional contact path geometries with vanishing W\"unschmann invariant as the spaces of null paths of a three-dimensional Lorentzian conformal structure, and realizing every three-dimensional Lorentzian conformal structure as the space of contact paths of a three-dimensional contact path geometry with vanishing W\"unschmann invariant. This goes back to W\"unschmann, Chern, and E. Cartan's \cite{Cartan-Spanish}. In \cite{Sato-Yoshikawa} H. Sato and A. Yoshikawa solve the equivalence problem for three-dimensional contact path geometries and explain Chern's results. One implication relevant here is that split quaternionic contact structures are a natural higher dimensional generalization of three-dimensional Lorentzian structures, something already apparent from the isomorphism $\smp(2, \rea) \simeq \mathfrak{so}(2, 3)$.

The author thanks Robin Graham for advice and suggestions. Andreas \v{C}ap suggested that it would be interesting to study the contact path geometries and made many useful comments; the author thanks him for this and for explaining his point of view. Some version of many of the results of Section \ref{prolongsection} appeared in the author's thesis, \cite{Fox-thesis}. 

\section{Contact Path Geometries}\label{prolongsection}
\subsection{Definitions and Notation}\label{contactpathbasics}

Let $M$ be a smooth manifold of dimension $n$, let $\pi:\ptm \to M$ be its projectivized tangent bundle, let $\piz:\tmz \to M$ be tangent bundle with the zero section deleted, and regard the defining projection, $\rho:\tmz \to \ptm$, as a principal $\rea^{\times}$ bundle. A \textbf{generalized conformal structure} (of \textbf{rank} $r$ and \textbf{degree} $1$) on $M$ is an open subset, $\C$, of a smooth, closed $(n+r)$-dimensional submanifold of $\ptm$ such that $\pi:\C \to M$ is a submersion. %, i.e. $\pi_{\ast}(L)(\C) = T_{\pi(L)}M$.
Note that by definition an open subset $\B \subset \C$ of a generalized conformal structure on $M$ is a generalized conformal structure over $\pi(\B) \subset M$. Let $E \subset T\C$ be the rank $r + 1$ tautological subbundle with fiber $E_{L} = \{X \in T_{L}\C: \pi_{\ast}(L)(X) \in L \subset T_{\pi(L)}M\}$. The vertical subbundle, $V = \ker \pi_{\ast} \subset T\C$, is a codimension one subbundle of $E$. An \textbf{automorphism} of the generalized conformal structure, $\C$, is a fiber preserving diffeomorphism of $\C$ which preserves $E$. The lift to $\ptm$ of a diffeomorphism of $M$ is the projectivization of its differential, and this lift preserves $E_{\ptm}$. Every diffeomorphism of $M$ the differential of which preserves $\C$ lifts to an automorphism of $\C$, %(an automorphism of $(M, \C)$) lifts to a fiber preserving diffeomorphism of $\C$ which preserves $E$,
and in many cases of interest it can be proved that every automorphism of $\C$ which preserves $E$ arises as such a lift. On the other hand, for a general choice of $\C$ there may be no diffeomorphisms of $M$ preserving $\C$.

\begin{definition}\label{nullpathdefn}
A \textbf{null path geometry} on $M$ is a foliation by one dimensional submanifolds transverse to the vertical of a generalized conformal structure $\C \subset \ptm$. Equivalently a null path geometry is a splitting $E = V \oplus W$ defined along $\C$. The projections into $M$ of the integral manifolds of $W$ are called the \textbf{null paths} of the null path geometry. Two null path geometries, $(\pi:\C\to M, W)$ and $(\pi^{\prime}:\C^{\prime}\to M^{\prime}, W^{\prime})$ are \textbf{equivalent} if there is a fiber preserving diffeomorphism from $\C$ to $\C^{\prime}$ preserving $E$ and mapping $W$ to $W^{\prime}$. A null path geometry the null paths of which are among the unparameterized geodesics of some affine connection is called a \textbf{null projective structure}. Two null projective structures are equivalent if they are equivalent as null path geometries. The null paths of a null projective structure are called \textbf{null geodesics}.
\end{definition}

\begin{remark}
Null path geometries are defined locally on the fibers as well as on the base. Often locality on the fibers will be suppressed in both the notation and the exposition. A null path geometry is said to be \textbf{supported} by the generalized conformal structure, $\C$, if it is defined on an open subset of $\C$.  The null path geometries supported by the generalized conformal structure $\C = \ptm$ are called simply \textbf{path geometries}. 
\end{remark}

\begin{remark}
Essentially everything to follow has a sensible (and immediate) analogue in the holomorphic category. In the holomorphic category there are interesting rigidity phenomena, e.g. a holomorphic path geometry defined on all of $\ptm$ must, by virtue of Hartogs' extension principle applied on the fibers of $\tmz$, be a projective connection.
\end{remark}

\begin{remark}
Degree $k$ null path geometries can be defined by replacing $\ptm$ by the Grassmann bundle $Gr(k, TM)$, defining $V$ and $E$ as above, and requiring additionally that $W$ be integrable. For related discussions see the first chapter of Y. Manin's \cite{Manin}, and various papers of S. Gindikin, e.g. the papers with J. Bernstein, \cite{Bernstein-Gindikin}.
\end{remark}

A null projective structure may be identified with the equivalence class, $[\nabla]$, comprising all affine connections on $M$ having among their unparameterized geodesics the given null paths. Given a generalized conformal structure, $\C$, an affine connection, $\nabla$, is said to \textbf{admit a full set of null paths} if every (unparameterized) geodesic of $\nabla$ tangent to $\C$ at one point is everywhere tangent to $\C$. If $\C$ is the projectivization of a rank $r$ distribution, $H \subset TM$, then $\nabla$ admits a full set of null paths if and only if for every choice of $\theta^{1}, \dots, \theta^{n-r}$ spanning the annihilator of $H$ there exist one-forms, $\beta^{i}_{j}$, such that $\sym \nabla \theta^{i} = \sym (\beta^{i}_{j}\tensor \theta^{j})$.

A null path geometry, $(\C, W)$, is \textbf{subordinate} to a null path geometry, $(\C^{\prime}, W^{\prime})$ if the null paths of $(\C, W)$ are among the null paths of $(\C^{\prime}, W^{\prime})$. In general it is easy to see that there are many inequivalent path geometries subordinate to a given null path geometry. The vertical bundle, $V$, and the tautological bundle, $E$, on $\C$ are the restrictions to $\C$ of the corresponding bundles on $\ptm$; simply extend the splitting $E = V \oplus W$ locally off of the submanifold $\C \subset \ptm$.

The null path geometries (resp. null projective structures) on the generalized conformal structure $\C = \ph$ over a $(2n-1)$-dimensional contact manifold, $(M, H)$, are called \textbf{contact path geometries} (resp. \textbf{contact projective structures}), and their null paths are called \textbf{contact paths}. Since a contactomorphism preserves $H$, its differential acts naturally on $\ph$. From the definition of $E$ it is easy to see that the bundle automorphism of $\ph$ induced by a contactomorphism of $(M, H)$ preserves $E$. Moreover, the lift to $\ph$ of a contact path in $M$ is necessarily tangent to $E$. Lemmas \ref{bracketgeneratingcontact} and \ref{liftlemma} below show that $E$ is $2$-step bracket generating in a precisely specified way, and that a fiber preserving diffeomorphism of $\proj(H)$ preserving $E$ is induced by a contactomorphism of $(M, H)$.

In addition to the bundles $V$ and $E$ already defined, define vector bundles $E^{\perp}$ and $\Hup = \pi_{\ast}^{-1}(H)$ with fibers
\begin{align*}
%&V_{L} = \{v \in T_{L}\proj(H): \pi_{\ast}(L)(v) = 0\},& &\rank V = 2n-3,&\\
%&E_{L} = \{v \in T_{L}\proj(H): \pi_{\ast}(L)(v) \in L \subset H_{\pi(L)}\},& &\rank E = 2n-2,&\\
&E_{L}^{\perp} = \{v \in T_{L}\proj(H): \pi_{\ast}(L)(v) \in L^{\perp} \subset H_{\pi(L)}\},& &\rank E^{\perp} = 4n-6,&\\
&\Hup_{L} = \{v \in T_{L}\proj(H)): \pi_{\ast}(L)(v) \in H_{\pi(L)}\},& &\rank \Hup = 4n-5.& 
\end{align*}
Here $L^{\perp}$ indicates the skew complement of $L$ in $H_{\pi(L)}$ with respect to the conformal symplectic structure on $H_{\pi(L)}$. 

The quotient $\Hup/V$ is isomorphic to the pullback bundle, $\pi^{\ast}(H)$, and so $\Hpull$ admits a tautological filtration $E/V \subset E^{\perp}/V \subset \Hpull$. The conformal symplectic structure on $H$ induces on $\Hpull$ a fiberwise conformal symplectic structure, and $E^{\perp}/V$ is the skew complement of $E/V$ in $\Hpull$. Canonically, $V_{L} = T_{L}(\proj(H_{\pi(L)})) \simeq \hom(L, H_{\pi(L)}/L)$, so that $V \simeq (\Hup/E)\tensor (E/V)^{\ast}$. A distinguished subbundle $U \subset V$ is defined as the inverse image of $(E^{\perp}/E)\tensor (E/V)^{\ast}$ under this identification. Precisely, the fiber $U_{L}$ is the subspace $\hom(L, L^{\perp}/L)$ of $T_{L}(\proj(H_{\pi(L)})) \simeq \hom(L, H_{\pi(L)}/L)$. Geometrically, each fiber, $\proj(H_{L})$, is the projectivization of a conformal symplectic vector space, so admits a canonical contact structure, which is $U_{L}$ (in this regard see section \ref{flatmodel}). The filtration 
\begin{equation}\label{canonicalcontactfiltration}
U \subset V \subset E \subset E^{\perp} \subset \Hup \subset T\proj(H),
\end{equation}
determined by these vector bundles depends only on the contact structure on $M$. If $2n-1 = 3$, then $E^{\perp} = E$ and $U$ is trivial.

The geometric meaning of the bundle $E^{\perp}$ is given by the observation that if $N \subset M$ is an isotropic submanifold in the sense that $TN \subset H_{|N}$, then $\proj(TN) \subset \ph$ is tangent to $E^{\perp}$. To prove this, let $L \in \proj(TN)$ and $v \in T_{L}\proj(TN)$. Because $L \subset T_{\pi(L)}N$ and $N$ is isotropic, there holds $T_{\pi(L)}N \subset T_{\pi(L)}N^{\perp} \subset L^{\perp}$, so that $\pi_{\ast}(L)(v) \in T_{\pi(L)}N \subset L^{\perp}_{\pi(L)}$, which shows that $\proj(TN)$ is tangent to $E^{\perp}$.

\begin{remark}
One can consider \textbf{contact $k$-path} geometries, which are defined as foliations of the isotropic Grassmann bundle $IG(k, H)$, or equivalently as splittings $E = V \oplus W$, (with $W$ integrable), where $E$ is the usual tautological bundle. As in the $k = 1$ case there is a filtration of $TIG(k, H)$ by subbundles $U \subset V \subset E \subset E^{\perp}$. The $k = n-1$ case is simpler because $E = E^{\perp}$. When $1 < k < n-1$ the story will be much like the $k = 1$ case, though the details remain to be worked out.
\end{remark}

\subsection{Flat Model for Contact Path Geometries}\label{flatmodel}
The projectivization, $\rho:\standrepz = \standrep - \{0\} \to \proj(\standrep)$ of a $2n$-dimensional real symplectic vector space, $(\standrep, \Omega)$, has a canonical contact structure defined by declaring the contact hyperplane at the one-dimensional subspace $L \in \proj(\standrep)$ to be the image under $\rho_{\ast}$ of the skew complement of $L$. The total space of the projectivization of the contact hyperplane is identified with the isotropic flag manifold of isotropic two-flags consisting of a one-dimensional subspace (the \textbf{pole}) contained in a two-dimensional isotropic subspace (the \textbf{cloth}). The leaves of the foliation of the isotropic flag manifold determining the contact path geometry are the submanifolds consisting of isotropic flags with a given fixed cloth. The symplectic group, $G = Sp(n, \rea)$, acts transitively on the set of isotropic subspaces of $(\standrep, \Omega)$ of a fixed rank $k \leq n$. The manifold of such subspaces is the isotropic Grassmannian $Q_{k} = IG(k, n, \rea) = G/P_{k}$, where $P_{k}$ is a parabolic subgroup that may be taken to be the stabilizer of the isotropic subspace spanned by $e_{1}, \dots, e_{k}$. It is convenient also to write $P_{i_{1}\dots i_{k}} = P_{i_{1}}\cap \dots \cap P_{i_{k}}$. 

$Q_{1}$ is simply $(2n-1)$-dimensional projective space, and $Q_{2}$ is the $(4n-5)$ dimensional isotropic Grassmannian of two-dimensional isotropic subspaces. The $(4n-4)$-dimensional space $G/P_{12} = Q_{12} = \{(V_{1}, V_{2}) \in Q_{1}\times Q_{2}: V_{1} \subset V_{2}\}$ is the space of isotropic $2$-flags $V_{1} \subset V_{2}$, where $V_{1}$ is $1$-dimensional and $V_{2}$ is $2$-dimensional and isotropic. 
%This gives the double fibration
%\begin{equation}\label{cpathdoublefibration}
%\xymatrix{
%&Q_{12} = \{(V_{1}, V_{2}) \in Q_{1}\times Q_{2}: V_{1} \subset V_{2}\} \ar[dl]^{\rho_{1}} \ar[dr]_{\rho_{2}}\\
%Q_{1}& & Q_{2}}
%\end{equation}
$G$ acts on each of the spaces $Q_{12}$, $Q_{1}$, $Q_{2}$, and the actions commute with the projections $\rho_{1}$ and $\rho_{2}$. The Maurer-Cartan form on $G$ determines on each of these homogeneous spaces a canonical $(\g, P)$ Cartan connection, where $P = P_{1}$, $P_{2}$, or $P_{12}$, as is appropriate. 

The fibers of the projections $\rho_{1}$ and $\rho_{2}$ form a pair of transverse foliations of $Q_{12}$. At the flag $V_{1} \subset V_{2}$ in $Q_{12}$ the leaves of the foliations may be described explicitly. The fiber of $\rho_{2}$ passing through $V_{1} \subset V_{2}$ comprises all isotropic $2$-flags with the cloth $V_{2}$, and is therefore identified with $\rpone$. The fiber of $\rho_{1}$ passing through $V_{1} \subset V_{2}$ comprises all isotropic $2$-flags with the pole $V_{1}$, and is therefore identified with $\proj(V_{2}^{\perp}/V_{2}) \simeq \mathbb{P}^{2n-3}(\rea)$. This shows that the total space of the $(4n-4)$-dimensional space $Q_{12}$ may be realized as the projectivization of the bundle of contact hyperplanes on $Q_{1}$. The images under $\rho_{1}$ of the fibers of $\rho_{2}$ form the $(4n-5)$-parameter family of contact paths in $Q_{1}$. The $2$-step bracket generating distribution, $E \subset TQ_{12}$, is defined by letting $E_{q}$, at a point $q \in Q_{12}$, equal the span of the tangent spaces to the fibers of $\rho_{1}$ and $\rho_{2}$ passing through the point $q$. The total space of the $(4n-4)$-dimensional space $Q_{12}$ is realized as the tautological $\rpone$-bundle over $Q_{2}$, the fiber over a $2$-dimensional isotropic subspace (a point in $Q_{2}$) comprising all $1$-dimensional subspaces contained in that subspace.  The images under $\rho_{2}$ of the fibers of $\rho_{1}$ form a $(2n-1)$-parameter family of copies of $\mathbb{P}^{2n-3}(\rea)$ in $Q_{2}$.

The flat model contact path geometry is the foliation of $Q_{12}$ by the fibers of $\rho_{2}$, and the Cartan connection to be associated to a contact path geometry is modeled on the canonical Cartan connection on $Q_{12} = G/P_{12}$. An element of $Q_{12}$ is an isotropic flag $V_{1} \subset V_{2}$. Suppose fixed $V_{1} \in Q_{1}$ and choose any $V_{1} \subset V_{2}$. The submanifold of $Q_{12}$ consisting of all flags with cloth $V_{2}$ is transverse to the fibers of $\rho_{1}$. Its image in $Q_{1}$ is a projective line, the space of one-dimensional subspaces of $V_{2}$. Each choice of $V_{2}$ containing $V_{1}$ determines such a submanifold of $Q_{1}$ passing through $V_{1}$. The union of the spaces tangent at $V_{1}$ to these submanifolds forms a linear subspace in $T_{V_{1}}Q_{1}$ which is the contact hyperplane on $Q_{1}$. This shows that, geometrically, the contact distribution $H$ at the one-dimensional subspace $V_{1}\subset \standrep$ is the image, in $Q_{1} = \proj(\standrep)$, of the $\Omega$-skew complement of $V_{1}$. There is a dual construction of a smoothly varying family of quadratic cones in the tangent spaces to $Q_{2}$. Suppose fixed $V_{2} \in Q_{2}$ and choose any $V_{1} \subset V_{2}$. The submanifold of $Q_{12}$ consisting of all isotropic $2$-flags with pole $V_{1}$ is transverse to the fibers of $\rho_{2}$. Its image in $Q_{2}$ is a submanifold containing $V_{2}$. Each choice of $V_{1}$ contained in $V_{2}$ determines such a submanifold of $Q_{2}$ passing through $V_{2}$. The union of the spaces tangent at $V_{2}$ to each of these submanifolds forms a cone in $T_{V_{2}}Q_{2}$. The possible choices of $V_{1}$ form a copy of the projective line. To each such choice is associated the tangent space to the submanifold of $Q_{2}$ determined by the choice of $V_{1}$. These tangent spaces rule the cone, and each may be taken as a generator of the cone. 

Tautological vector bundles, $L$ and $L^{\perp}$, over $Q_{k}$ are defined by
\begin{align*}
&L = \{(L, v) \in Q_{k} \times \standrep: v \in L\},& &L^{\perp} = \{(L, v) \in Q_{k} \times \standrep: v \in L^{\perp}\}.& 
\end{align*}
The symplectic form $\Omega$ induces an isomorphism $\standrep/L \simeq (L^{\perp})^{\ast}$. Using the isotropy condition it is straightforward to see that tangent bundle $TQ_{k}$ is canonically identified with the subspace of $\hom(L, \standrep/L) \simeq \standrep/L \otimes L^{\ast} \simeq (L^{\perp})^{\ast}\tensor L^{\ast}$ consisting of elements the restrictions of which to $L \tensor L$ are in the subspace $S^{2}(L^{\ast}) \subset (L^{\perp})^{\ast}\tensor L^{\ast}$. If $k < n$ then $L$ is contained strictly in $L^{\perp}$, and there is a non-trivial subbundle of $TQ_{k}$ defined by $C = (L^{\perp}/L)\tensor L^{\ast} \simeq \hom(L, L^{\perp}/L)$. The isomorphism $\standrep/L \simeq (L^{\perp})^{\ast}$ induced by $\Omega$ identifies the subbundle $L^{\perp}/L$ with the subspace of $(L^{\perp})^{\ast}$ annihilating $L$, so that $C$ is identified with the subspace of $(L^{\perp})^{\ast}\tensor L^{\ast}$ the elements of which vanish when restricted to $L \tensor L$. Restriction of maps induces a projection $(L^{\perp})^{\ast} \to L^{\ast}$, and $C$ is the kernel of the resulting projection $TQ \to S^{2}(L^{\ast})$, which shows also that $TQ_{k}/C \simeq S^{2}(L^{\ast})$. It follows that $C$ has rank $2k(n-k)$ and corank $\frac{k(k+1)}{2}$, so that the dimension of $Q_{k}$ is $\frac{k(k+1)}{2} + 2k(n-k)$. In the case $k = n$, $L^{\perp} = L$, and so $TQ_{n} \simeq S^{2}(L^{\ast})$, i.e. the tangent bundle of the Lagrangian Grassmannian, $Q_{n}$, is identified with the symmetric $2$-forms on $L$ (when $n = 2$, this is the case relevant for the $3$-dimensional contact path geometries). When $k = 1$, $C$ is a contact structure, the preceeding discussion recovers the usual identification of the bundle of positive contact one-forms with the square of the bundle of frames in the tautological bundle, $L$. %SHOWN EXPLICITLY LATER

There is also some distinguished algebraic structure in $TQ_{k}$. Consider, for each $0 \leq l < k$, the cones in $\hom(L, V/L)$ and in $\hom(L, L^{\perp}/L)$ consisting of maps with rank less than or equal to $l$. These are smoothly varying cones of degree $l$ (in coordinates, they are described by the vanishing of all $(l+1) \times (l+1)$ minors of the corresponding linear transformation), which are non-trivial when $k \neq 1, n$.

The bundle $L^{\perp}/L$ has a symplectic structure, $\omega$, defined by $\omega(u+L, v+L) = \Omega(u, v)$, and this induces on $C$ an $S^{2}(L^{\ast})$-valued symplectic structure defined by $\omega((u+L)\tensor a, (v + L)\tensor b) = \Omega(u, v)a \odot b$ (which locally can be regarded as a symmetric-matrix valued conformal symplectic structure). For $1 \leq k \leq n-1$, the fiber over $L \in Q_{k}$ of the projection $Q_{kl} \to Q_{k}$ is $IG(k, L^{\perp}/L) \simeq IG(k, n-k, \rea)$. The space $IG(k, L^{\perp}/L)$ has itself a canonical multicontact structure, and these multicontact structures on the fibers $IG(k, L^{\perp}/L)$ determine on $Q_{kl}$ a subbundle, $U$, of the vertical bundle of the projection $Q_{kl} \to Q_{l}$. Moreover, from the description of the multicontact structure, $C$, it is apparent that the brackets of sections of $U$ generate the vertical. This subbundle $U$ is a non-obvious feature of the geometry of $Q_{kl}$ which plays an important role in the context of curved contact path geometries.

Associated to any distribution, $C$, on the smooth manifold $N$, there are a $TN/C$-valued one form, $\btheta$, and a section, $\bomega$, of $\Lambda^{2}(C^{\ast}) \tensor \ann(C)^{\ast}$ defined as follows. There is a canonical isomorphism $TN/C \simeq \ann(C)^{\ast}$ defined by $v + C \to \alpha_{v}$ where, for $\theta \in \ann(C)$, $\alpha_{v}(\theta) = \theta(v)$. The inclusion $\ann(C) \to \Lambda^{1}(T^{\ast}N)$ is linear so may be viewed as a section, $\btheta$, of $\Lambda^{1}(T^{\ast}N)\tensor \ann(C)^{\ast}$. Using the canonical isomorphism $TN/C \simeq \ann(C)^{\ast}$, $\btheta$ may be viewed as a section of $\Lambda^{1}(T^{\ast}N)\tensor TN/C$. The map $\ann(C) \to \Lambda^{2}(C^{\ast})$ defined by $\theta \to d\theta|_{C \times C}$ is linear, so may be viewed as a section, $\bomega$, of $\Lambda^{2}(C^{\ast}) \tensor \ann(C)^{\ast}$. Using the canonical isomorphism $TN/C \simeq \ann(C)^{\ast}$, $\bomega$ may be viewed as a section of $\Lambda^{2}(C^{\ast})\tensor TN/C$. Composition of the canonical projection $\rho:TN \to TN/C$ with the Lie bracket $C \times C \to TN$ induces a tensorial mapping $\levi:\Lambda^{2}(C) \to TN/C$. Evidently, $-\levi = \bomega$. 

Let $C$ be the canonical multicontact structure on the isotropic Grassmannian, $Q = Q_{k}$. Using the canonical isomorphism $TQ/C \simeq S^{2}(L^{\ast})$, $\btheta$ and $\bomega$ may be viewed as taking values in $S^{2}(L^{\ast})$ (so that $\bomega$ is the $S^{2}(L^{\ast})$-valued symplectic structure on $C$ constructed above). A trivialization of $TQ/C$ allows $\btheta$ and $\bomega$ to be represented locally as forms taking values in symmetric $k \times k$ matrices, as will be described now. Choose coordinates $x^{\alpha}, x_{-\alpha}, x^{p}$ on $\standrep$ and write the symplectic structure on $\standrep$ as $\Omega = dx^{\alpha}\wedge dx_{-\alpha} + \tfrac{1}{2}\omega_{pq}dx^{p}\wedge dx^{q}$, where lowercase Latin indices run over $\{1, \dots, 2n-2k\}$, lowercase Greek indices run over $\{1, \dots, k\}$, and capital Latin indices run over $\{1, \dots, 2n\}$. In these local coordinates, a path, $\gamma(t) \in Q_{k} = IG(k, n, \rea)$, may be written as the $2n \times k$ matrix
\begin{align*}
\gamma = \begin{pmatrix} \delta_{\al}\,^{\be} \\ x_{\alpha}\,^{p} \\ z_{\alpha\beta}\end{pmatrix},
\end{align*}
The condition that such a path lie in $Q_{k}$ is $\Omega(\gamma, \gamma) = 0$, which is equivalent to the equations
\begin{align*}
z_{\al\beta} - z_{\be\alpha} + \omega_{pq}x_{\alpha}\,^{p}x_{\beta}\,^{q} = 0.
\end{align*}
%OMITTED TO SAVE SPACE
%Differentiating along $\gamma$ gives
%\begin{align*}
%\dot{z}_{\al\beta}  + \omega_{pq}x_{\alpha}\,^{p}\dot{x}_{\beta}\,^{q}= \dot{z}_{\be\alpha} + \omega_{pq}x_{\beta}\,^{p}\dot{x}_{\alpha}\,^{q}.
%\end{align*}
The distribution $C$ is defined by the condition that the instantaneous motion at a point be in the skew-complement of the point, i.e. $\Omega(\dot{\gamma},\gamma) = 0$. In coordinates this is expressible as
\begin{align*}
&\dot{z}_{\al\beta} + \omega_{pq}x_{\alpha}\,^{p}\dot{x}_{\beta}\,^{q} = 0 = \dot{z}_{\be\alpha} + \omega_{pq}x_{\beta}\,^{p}\dot{x}_{\alpha}\,^{q}.
\end{align*}
This shows that $C$ is the kernel of the $S^{2}(L^{\ast})$-valued one-form, $\Theta = \theta_{\al\beta}$, expressible locally as
\begin{align*}
&\theta_{\al\beta} = dy_{\al\be} + \tfrac{1}{2}\omega_{pq}x_{\al}\,^{p}dx_{\be}\,^{q} + \tfrac{1}{2}\omega_{pq}x_{\be}\,^{p}dx_{\al}\,^{q},
\end{align*}
where the coordinates $y_{\al\be} = y_{\be\al}$ are defined by $y_{\al\be} = z_{\al\be} - \tfrac{1}{2}\omega_{pq}x_{\al}\,^{p}x_{\be}\,^{q}$. By the kernel of $\Theta$ is meant $\cap_{\alpha, \beta = 1}^{k}\ker\theta_{\alpha\be}$. Evidently, $\Theta$ is the specialization of $\btheta$ determined by a local trivialization of $S^{2}(L^{\ast})$ and so the matrix $\Theta_{\al\be}$ is determined only up to similarities of the form $c_{\al}\,^{\ga}c_{\be}\,^{\si}\Theta_{\ga\si}$ where $c_{\al}\,^{\be}$ is a smooth $GL(k, \rea)$-valued function on $Q_{k}$ (the $c_{\al}\,^{\be}$ should be thought of as the local expressions of a $GL(k, \rea)$ principal connection on the bundle of frames in $L$).

The vector fields $X_{i}\,^{\al} = \tfrac{\pr}{\pr x_{\al}\,^{i}} + \omega_{ip}x_{\sigma}\,^{p}\tfrac{\pr}{\pr y_{\sigma \al}}$ span $C$ and satisfy the bracket relations $[X_{i}\,^{\al}, X_{j}\,^{\be}] = -2\omega_{ij}V^{\al\be}$ where $V^{\al\be} = \tfrac{\pr}{\pr y_{\al\be}}$. All the other Lie brackets vanish. In particular, this shows that $C$ is $1$-step bracket generating. These bracket relations are encoded also in the matrix valued $2$-form $\Omega_{\al\be} = d\theta_{\al\be} = \omega_{ij}dx_{\al}\,^{i}\wedge dx_{\be}\,^{j}$, which is the specialization of $\bomega$ determined by a local trivialization of $S^{2}(L^{\ast})$. While $\Omega_{\al\be}$ depends on the choice of $\Theta_{\al\be}$, its restriction to $C \times C$ is well-defined up to similarities of the form $c_{\al}\,^{\ga}c_{\be}\,^{\si}\Omega_{\ga\si}$ where $c_{\al}\,^{\be}$ is a smooth $GL(k, \rea)$-valued function on $Q_{k}$.

\begin{remark}
%%REDUCED FOR SPACE
%The $(4n-5)$-dimensional manifold, $J^{1}(\rea^{n-2}, \rea^{3})$ of $1$-jets of smooth functions from $\rea^{n-2}$ to $\rea^{3}$ admits a tautological codimension $3$ bracket-generating distribution given in natural coordinates, $\{(z^{\al}, s^{i}, p^{\al}_{i}): 1 \leq \al \leq 3$, $1 \leq i \leq n-2\}$, as the annihilator of the span of the $1$-forms $\phi_{\al} = dz^{\al} - p^{\al}_{i}ds^{i}$. This distribution is not smoothly equivalent to the distribution $C$ on $Q_{2}$ as may be shown by computing the wedge length of its annihilator (see \cite{Bryantetal}). Precisely, $(d\phi_{\al})^{n-2} \neq 0$ and $(d\phi_{\al})^{n-1} = 0$. On the other hand, for $\al \neq \be$ the one-forms $d\theta_{\al\be}$ on $Q_{2}$ satisfy that $(d\theta_{\al\al})^{n-2} \wedge (d\theta_{\be\be})^{n-2}$ is a non-zero multiple of $(d\theta_{\al\be})^{2n-4} \neq 0$. 
The contact structure on the space of $1$-jets of smooth functions from $\rea^{n-2}$ to $\rea^{3}$ is not locally smoothly equivalent to $C$ on $Q_{2}$ as may be shown by computing the wedge lengths of their annihilators (see \cite{Bryantetal}).
\end{remark}

A $\Det(S^{2}(L^{\ast}))$-valued $2k$-form, $\bpsi$, on $C$ is defined by taking the determinant of the $S^{2}(L^{\ast})$-valued $2$-form, $\bomega$. $\bpsi$ may be regarded as a conformal class, $[\Psi]$, of $2k$-forms on $C$ by letting $\Psi$ equal the restriction to $C$ of $\det_{\eta} \Omega_{\al\be}$, where $\eta \in \Lambda^{k}(L)$. In coordinates the conformal class $[\Psi]$ is represented by 
\begin{align*}
\Psi = \eta^{\al_{1}\dots \al_{k}}\eta^{\be_{1}\dots \be_{k}}\omega_{i_{1}j_{1}}\dots\omega_{i_{k}j_{k}}dx_{\al_{1}}\,^{i_{1}} \wedge dx_{\be_{1}}\,^{j_{1}}\wedge \dots \wedge dx_{\al_{k}}\,^{i_{k}} \wedge dx_{\be_{k}}\,^{j_{k}}.
\end{align*}
It can be checked that $\bpsi$ is non-degenerate in the sense that the restriction to $C$ of the $2k(n-k)$-form $\Psi^{n-k}$ is non-vanishing.

The set of $k$-dimensional subspaces of an $n$-dimensional isotropic subspace of $\standrep$ is a copy of $G(k, n, \rea)$ embedded in $Q_{k}$. Its dimension is $k(n-k)$ and it is tangent to the multicontact structure $C$ on $Q_{k}$. If for $1 \leq I \leq n-k$ there are chosen constants $C_{I}\,^{i}$ so that $C_{I}\,^{i}$ has rank $n-k$ and $C_{I}\,^{i}C_{J}\,^{j}\omega_{ij} = 0$, then it is easily checked that $A_{I}\,^{\al} = C_{I}\,^{i}X_{i}\,^{\al}$ are $k(n-k)$ vector fields spanning an integrable $[\Psi]$-null subbundle of $C$. The integral manifolds are the copies of $G(k, n, \rea)$ just described.

\subsubsection{Special Structures on $Q_{2}$}
The principal special feature of $Q_{2}$ is that the subbundle, $C \subset TQ_{2}$, supports a canonical a conformal class, $[g]$, of metrics of split signature. 

Let $(\standrep, \Omega)$ be a symplectic vector space, let $L$ be a $2$-dimensional vector space, and set $\standw = \standrep \tensor L^{\ast}$. Denote by $\eta$ any non-zero section of $\Lambda^{2}(L)$. Define $g \in \Gamma(S^{2}(\standw^{\ast}))$ by $g(x \tensor a, y \tensor b) = \Omega(x, y)\eta(a, b)$. Varying the choice of $\eta$ rescales $g$, so that the conformal class, $[g]$, is well-defined. Non-degeneracy of $[g]$ follows from non-degeneracy of $\eta$ and $\Omega$. If $J \subset \standrep$ is $\Omega$-isotropic then $J \tensor L^{\ast}$ is $[g]$-null, which shows that $[g]$ admits null planes of dimension $\dim \standrep$, so that $[g]$ has split signature. Applying this construction to $C \subset TQ_{2}$ using the identification $C \simeq L^{\perp}/L \tensor L^{\ast}$ and the fiberwise symplectic structure induced by $\Omega$ on $L^{\perp}/L$, yields a conformal class, $[g]$, of signature $(2n-4, 2n-4)$ metrics on $C$. 

%ONMITTED TO SAVE SPACE
%A central simple associative algebra is called quaternionic if it has dimension $4$. A quaternionic algebra is called split if it is isomorphic to the algebra of matrices over the base field. Because a finite-dimensional central simple algebra is simple artinian, a quaternionic algebra is by the Artin-Wedderburn theorems isomorphic to a matrix ring over a division algebra over the base field, and so if the base field is $\rea$ it is easy to see that a quaternionic algebra is split if and only if it is isomorphic to the endormorphism algebra of a two-dimensional real vector space, $\standu$. 

The \textbf{split quaternion algebra}, $\squat$, is the endomorphism algebra of a two-dimensional real vector space, $\standu$. Each $\eta \in \Lambda^{2}(\standu)$ induces an isomorphism $\standu \to \standu^{\ast}$. Conjugating the transpose of $a \in \squat$ by this isomorphism gives an involution of $\squat$ denoted by $a \to \bar{a}$. With respect to a unimodular basis of $(\standu, \eta)$, $\bar{a}$ is the matrix of cofactors so that a norm is defined on $\squat$ by $|a|^{2}I = \bar{a}a = \det_{\eta}(a)I $, where $I$ is the unit in $\squat$. The associated quadratic form $G(a, b) = \re \bar{a}b$ deteremines on $\squat$ a signature $(2, 2)$ metric the conformal class of which depends only on the algebra structure of $\squat$, and not on the choice of $\eta$. The real part of $a \in \squat$ is $\re a = \tfrac{1}{2}(a + \bar{a}) = \tfrac{1}{2}(\tr a) I$, and the imaginary part is $\im a = \tfrac{1}{2}(a - \bar{a}) = \tfrac{1}{2}(a - (\tr a)I)$. The imaginary elements, $\im \squat$, of $\squat$ are just the trace free elements. 

The group of elements of norm $1$ (resp. units) in $(\squat, \eta)$ is $SL(2, \rea)$ (resp. $GL(2, \rea)$) which can be regarded also as $Sp(1, \rea)$ (resp. $CSp(1, \rea)$), as these are just the elements which, when viewed as endomorphisms of $\standu$, preserve $\eta$. (Recall that the group of unit norm elements of the usual quaternions is the compact symplectic group $S^{3}$; the usual quaternions and the split quaternions are different real forms of $\mathbb{M}_{2}(\com)$, and their groups of unit norm elements are, respectively, the compact and the split real forms of $Sp(1, \com)$). The Lie bracket $[a, b] = ab - ba = 2\im(ab)$ identifies $\im \squat$ with $\smp(1, \rea) \simeq \sltwor$. The restriction of $G$ to $\im \squat$ is a metric of signature $(1, 2)$ which is preserved by the adjoint action of $Sp(1, \rea)$ on $\im \squat$. This identifies the adjoint group $PSp(1, \rea)$ with $SO(1, 2)$.

Because any two symplectic vector spaces of the same dimension are isomorphic, for any choice of $\eta$ there is a symplectic linear isomorphism $(\standrep, \Omega) \simeq (\oplus_{i = 1}^{n}L, \oplus_{i = 1}^{n}\eta)$, and this gives an identification of right $\squat$-modules $\standrep \tensor L^{\ast} \simeq (\oplus_{i = 1}^{n} L) \tensor L^{\ast} \simeq \oplus_{i = 1}^{n}(L \tensor L^{\ast}) = \squat^{n}$. Under the identification $\squat = L \tensor L^{\ast}$, a representative, $G \in [G]$, of the conformal metric on $\squat$ has the form $G(x \tensor a, y \tensor b) = \eta(x, y)\eta(a, b)$. This shows that the conformal class $[\oplus_{i = 1}^{n} G]$ is identified with the conformal class of metrics, $[g]$, on $\standrep \tensor L^{\ast}$ defined above. On the other hand, $\oplus_{i = 1}^{n} G$ is just the induced inner product $\re(\bar{p}q)$ for $p, q \in \squat^{n}$. 

As a right $\squat$-module automorphism is an $\rea$-linear map, it must preserve the tensor product decomposition $\standrep \tensor L^{\ast}$, and so must have the form $p \to (A, a)\cdot p = Ap\bar{a}$ where $A \in GL(\standrep)$ and $a \in \squat$ is invertible. Since $(A, a)$ must commute with the right action of $\squat$ on $\squat^{n}$, the element $a$ must be in the center of $\squat$, which is $\rea^{\times}$. Since for $a \in \rea^{\times}$, $(A, a)$ and $(a^{-1}A, 1)$ act in the same way on $\squat^{n}$ this shows that the right $\squat$-module automorphisms of $\squat^{n}$ are isomorphic to $GL(n, \rea)$. If such an automorphism preserves the conformal metric $[g]$ on $\squat^{n}$ then it follows from the discussion of the preceeding paragraph that $A$ preserves $\Omega$ up to a scalar factor, so that the group of $[g]$-conformal right $\squat$-module automorphisms of $\squat^{n}$ is isomorphic to the real conformal symplectic group $CSp(n, \rea)$. Using this identification the $\rea$-linear action on $\squat^{n}$ of the group $\cssp$ of equivalence classes of pairs $(A, a) \in CSp(n, \rea) \times \squat^{\times}$ is defined by $(A, a)\cdot p = Ap\bar{a}$ (the equivalence is $(A, a) \simeq (rA, r^{-1}a)$ for $r \in \rea^{\times}$). The group $CSp(n, \rea)$ acts on $\squat^{n}$ as conformal right $\squat$-module automorphisms, and the group $GL(2, \rea)$ acts conformally on $\squat^{n}$ on the right as the group of units in $\squat$; the intersection in $CO(\squat^{n})$ of these two groups is just $\rea^{\times}$, and the subgroup they generate is $\cssp$. The connected component of the identity $CSp_{0}(n, \rea)$ acts on $\squat^{n}$ as positive $[g]$-conformal right $\squat$-module automorphisms.
 
Every finite-dimensional right $\squat$-module is isomorphic to a direct sum of copies of $L$, so has even dimension. The rank of a non-zero element $x \in \standrep \tensor L^{\ast}$ is the dimension of the image of the $\rea$-linear map $x:L \to \standrep$. An element of $\squat^{n}$ is decomposable in the sense that it has the form $a \tensor b$ if and only if it has rank $1$. For $x \in \squat^{n}$, the right $\squat$-submodule $x \squat \subset \squat^{n}$ is identified with all linear maps $x \circ a:L \to \standrep$, where $a \in \squat$ is viewed as an endomorphism of $L$. Since all these linear maps have the same image, this shows that the right $\squat$-submodule generated by $x$ is naturally identified with a subspace of the form $K \tensor L^{\ast}$, where $K$ is the image of $x$ viewed as a linear map (note that $K$ may have dimension $1$ or $2$). $K$ has dimension $1$ if and only if $x$ is decomposable; in the event $K$ has dimension $2$ there may be chosen two rank $1$ elements of $\squat^{n}$ generating $x\squat$. Similarly, any $2k$-dimensional right $\squat$-submodule of $\squat^{n}$ is generated by $k$ rank $1$ elements of $\squat^{n}$, and so has the form $K \tensor L^{\ast}$ for the $k$-dimensional subspace $K \subset \standrep$ spanned by the images in $\standrep$ of these generating elements. The projective space $Q_{1} = \proj(\standrep)$ is identified with the space of $2$-dimensional right $\squat$-submodules of $\squat^{n}$. More generally, the ordinary Grassmannian, $Gr(k, \standrep)$ is identified with the space of $2k$-dimensional right $\squat$-submodules of $\squat^{n}$.

$\squat$ is generated over $\rea$ by $1, j, e, f$ subject to the relations $-j^{2} = e^{2} = f^{2} = 1$, and $ef = j$ (which imply $fe  = -j$). The standard representation of $\squat$ as $2 \times 2$ matrices is
\begin{align*}
\rho(a + bj + ce + df) = \begin{pmatrix}  a+ d & c - b\\ b + c & a - d\end{pmatrix}
\end{align*}
For example, the element $\tfrac{1}{2}(1 + f)$ is idempotent, and using this it is easily checked that the right $\squat$-module $(1 + f)\squat$ is $2$-dimensional, whereas $j\squat = \squat$. 

For $1 \leq k \leq n$ define $Gr(k, n, \squat)$ to be the space of $2k$-dimensional right $\squat$-submodules of $\squat^{n}$. The preceeding shows that $Gr(k, n, \squat)$ comprises the $2k$-dimensional subspaces of $\squat^{n}$ of the form $K \tensor L^{\ast}$, so that as smooth manifolds $Gr(k, n, \squat) \simeq Gr(k, \standrep)$. The $2k$-dimensional subspace $K \tensor L^{\ast}$ is $[g]$-null if and only if $K$ is isotropic. Denote by $Q^{k-s}_{k}$ the subset of $Gr(k, n, \squat)$ comprising submodules for which the the dimension of a maximal isotropic subspace of the associated subspace $K \in \standrep$ is $s \in \{k, k-2, \dots \}$. The disjoint union $Gr(k, n, \squat) = \cup_{s \leq k}Q^{k-s}_{k}$ is a partition into orbits of $Sp(n, \rea)$ (which are smooth submanifolds). The unique closed orbit is the isotropic Grassmannian $Q^{0}_{k} = Q_{k}$, and the dimension of $Q^{k-s}_{k}$ is $\binom{s+1}{2} + 2s(n-s) + (k-s)(2n -k -s)$. The maximal orbit (having $s = 0$ or $s = 1$ depending on the parity of $k$) is an open domain in $Gr(k, \standrep)$. There must be some sense in which $Q_{k}$ can be viewed as a sort of \v{S}ilov boundary of this maximal orbit.

Every two-dimensional subspace of $(\standrep, \Omega)$ is either isotropic or symplectic, so the symplectic group $Sp(n, \rea)$ has two orbits on $Gr(2, \standrep) = Gr(1, n, \squat)$. The closed orbit is $Q_{2}$, comprising the $g$-null points (isotropic subspaces) in $ Gr(1, n, \squat)$, and the other is its complement, which is identified with the $4(n-1)$-dimensional homogeneous space $\sgr \simeq Sp(n, \rea)/(Sp(n-1, \rea) \times Sp(1, \rea))$ comprising symplectic subspaces. This realizes $Q_{2}$ as a real codimension $1$ submanifold of $Gr(1, n, \squat)$ which should be regarded as the boundary at infinity of $\sgr$. For $u \tensor a, v \tensor b \in \standrep/\standu \tensor \standu^{\ast} \simeq T_{\standu}Gr(2, \standrep)$, setting $g(u\tensor a, v \tensor b) = \Omega(u, v)\Omega(a, b)$ defines on $Gr(2, \standrep)$ a symmetric two tensor which restricts to a split signature metric on $SGr(2, \standrep)$ (and vanishes along $IG(2, \standrep)$) such that $Sp(n, \rea)$ acts on $SGr(2, \standrep)$ by isometries.

An element $p \in \squat^{n}$ of rank $1$ has norm $0$. An element $p \in \squat^{n}$ of rank $2$ is either isotropic or symplectic; the isotropic rank $2$ elements are exactly the rank $2$ elements of norm $0$. To see the last claim, observe that because $L$ is of dimension $2$ it suffices to consider elements of the form $a \tensor u + b \tensor v$ with $u$ and $v$ linearly independent, and that the norm of such an element is a non-zero multiple of $\Omega(a, b)$. The symplectic Grassmannian, $SGr(2, \standrep)$, is the quotient by $SL(2, \rea)$ of the subspace of $\squat^{n}$ comprising unit norm elements. 

The tangent space to $Q_{2}$ at $p\squat$ is identified with the subspace of the vector space of $\squat$-module maps from $p\squat$ to $\squat^{n}/p\squat$ comprising those maps with image $[g]$-orthogonal to $p$. The fiber at $p\squat$ of the subbundle $C \subset TQ_{2}$ is identified with $\squat$-module maps having range in $(p\squat)^{\perp}/p\squat$, where $(p\squat)^{\perp}$ is the largest $\squat$-submodule of $\squat^{n}$ $[g]$-orthogonal to $p\squat$. Because $p\squat$ is $[g]$-null, $[g]$ descends to give a conformal metric on $C$. The graded bundle $C \oplus TQ_{2}/C$ is regular with fibers isomorphic as graded Lie algebras to $\g_{-}$; this Lie algebra is called the \textbf{split quaternionic Heisenberg Lie algebra}.

When $k = 2$ the conformal class of split signature metrics, $[g]$, on $C$ is expressible in coordinates as $g = \eta^{\al\be}\omega_{ij}dx_{\al}\,^{i} \tensor dx_{\be}\,^{j}$, for $\eta \in \Lambda^{2}(L^{\ast})$, so that $g(X_{i}\,^{\al}, X_{j}\,^{\be}) = \eta^{\al\be}\omega_{ij}$. Define $E, F, J \in \eno(C)$ satisfying the relations $E^{2} = F^{2} = I = -J^{2}$ and $E \circ F = J$ by
\begin{align*}
&E(X_{i}\,^{1}) = X_{i}\,^{2},& 
&F(X_{i}\,^{1}) = X_{i}\,^{1},& 
&J(X_{i}\,^{1}) = X_{i}\,^{2},& \\
& E(X_{i}\,^{2}) = X_{i}\,^{1},& & F(X_{i}\,^{2}) = -X_{i}\,^{2},& & J(X_{i}\,^{2}) = -X_{i}\,^{1},
\end{align*}
and such that when restricted to $C$ there hold $g(E(-), E(-)) = -g(-, -) = g(F(-), F(-))$ and $g(J(-), J(-)) = g(-, -)$. 
%$E$ and $F$ are reflections and $J$ is an almost complex structure. 
It is straightforward to check 
\begin{align*}
&\Omega_{12} = g(F\cdot, \cdot),& &\Omega_{11} = -g(E\cdot, \cdot) - g(J\cdot, \cdot),& &\Omega_{22} = g(E\cdot, \cdot) - g(J\cdot, \cdot),
\end{align*}
so that if $\eta$ is chosen so that $\eta^{12} = 1$ then the $\sltwor$-valued $2$-form $\eta^{\al\ga}\Omega_{\be\ga}$ has the form
\begin{align*}
&\Omega_{\be}\,^{\al} = \eta^{\al\ga}\Omega_{\be\ga} = \begin{pmatrix} \Omega_{1}\,^{1} & \Omega_{2}\,^{1}\\ \Omega_{1}\,^{2} & \Omega_{2}\,^{2}\end{pmatrix} = \begin{pmatrix}g(F\cdot, \cdot) & g(E\cdot, \cdot) - g(J\cdot, \cdot)\\ g(E\cdot,\cdot) + g(J\cdot,\cdot) & -g(F\cdot, \cdot)  \end{pmatrix}.
\end{align*}
An action of $\squat$ on $C$ is defined by associating to $a + bj + ce + df$ the endomorphism $aI + bJ + cE + dF$ of $C$. Regard $\squat^{n}$ as a right $\squat$-module. An imaginary idempotent in $\squat$ acts on $\squat^{n}$ as a \textbf{reflection} with respect to the inner product induced on $\squat^{n}$ by $G$. The elements of $\im \squat$ of norm act on $\squat^{n}$ as \textbf{complex structures}. To each reflection (resp. complex structure) $p \in \im \squat$ is associated a reflection (resp. almost complex structure) on $C$. This data provides the model for the split quaternionic contact structures discussed in Section \ref{splitquatstructure}. As a submanifold of $\im \squat$ the reflections (resp. complex structures) constitute a hyperboloid of one sheet (resp. two sheets), and observation which should be useful in attaching to a split quaternionic contact structure a twistor space.
%UNUSED SO OMITTED
%It is easily checked that the product of two orthogonal reflections is a complex structure. 

\subsection{Local Description of General Contact Path Geometries}\label{localdescriptioncontactpath}
Some notation, to be used whenever working in coordinates on contact path geometries, is fixed. Let capital Latin indices run over the range $0, 1, \dots, 2n-4, \infty$, lowercase Greek indices run over the range $0, \dots, 2n-4$ and lowercase Latin indices run over the range $1, \dots, 2n-4$. (When $2n-1 = 3$ interpret all expressions involving lowercase Latin indices as null). Let the constants $\omega_{ij}$ be skew-symmetric and non-degenerate, and let $\omega^{ij}$ satisfy $\omega^{ip}\omega_{pj} = -\delta_{j}\,^{i}$ and raise and lower lowercase Latin indices using $\omega_{ij}$, e.g. $x_{i} = x^{p}\omega_{pi}$. 

 If $A$ is abelian group the vector space, $V$, is \textbf{$A$-graded} if there is a direct sum decomposition $V = \bigoplus_{I \in A}V_{I}$. An $A$-grading on a Lie algebra, $\g$, is \textbf{compatible} with the Lie bracket if $[\g_{I}, \g_{J}] \subset \g_{I + J}$ for all $I, J \in A$. Each homomorphism of abelian groups, $\phi:A \to B$, associates to an $A$-graded vector space a $B$-graded vector space $V = \oplus_{b \in B}V_{b}$ defined by $V_{b} = \oplus_{I \in A: \phi(I) = b}V_{I}$. In this way the homomorphism $|\,|:\zk \to \integer$ defined by $|I| = \sum_{i_{p} \in I}i_{p}$ associates to each $\zk$-grading of $V$ an underlying $\integer$-grading. In this way each compatible $\zk$-grading of a Lie algebra $\g$ determines a compatible $\integer$-grading of $\g$.  For $I, J \in \zk$, write $I \leq J$ if for every $p \in \{1, \dots, k\}$, $i_{p} \leq j_{p}$, and write $I < J$ if, moreover, for some such $p$, there holds the strict inequality $i_{p} < j_{p}$. Write $\f_{I} = \oplus_{I \leq J}\g_{J}$, where $I \leq J$ refers to the partial order on $\zk$. 

If $\g$ is a semisimple (real or complex) Lie algebra and $\p \subset \g$ a subalgebra, the Killing form of $\g$ determines an orthogonal decomposition $\g = \g_{-}\oplus \p$. The subalgebra, $\p$, is called parabolic if it contains a maximal solvable subalgebra of $\g$, in which case there is a Killing orthogonal decomposition, $\p = \g_{0} \oplus \pplus$, where $\g_{0}$ is reductive and $\p^{+}$ is a nilpotent ideal identified via the Killing form with $\g_{-}^{\ast}$. There exists in $\g_{0}$ a unique \textbf{scaling} element, $E_{\p} \in \g_{0}$, characterized by the requirement that the eigenvalues of $\ad(E_{\p})$ acting on $\g_{-}$ be negative integers inclusive of $-1$ (so $\ad(E_{\p})$ has on $\p^{+}$ positive integer eigenvalues). The action of $E$ induces on any $\g_{0}$-module a $\g_{0}$-invariant $\integer$-grading, and if the $\g_{0}$ action extends to make the module a $\p$-module, there is on the module a $\p$-invariant $\integer$-filtration. From the geometric point of view the structures of signficance are these $\integer$-filtrations, but from the algebraic point of view it is often convenient to work with the associated graded objects. The nilpotent subalgebra $\g_{-}$ is a $\g_{0}$-module, and the decomposition of $\g_{-}$ into $\ad(E_{\p})$ eigenspaces determines a $\integer$-grading, $\g_{-} = \oplus_{i = -1}^{-k}\g_{i}$, and a Killing dual $\integer$-grading, $\p^{+} = \oplus_{i =1}^{k}\g_{i}$. The Lie algebra $\g = \oplus_{i = -k}^{k}\g_{i}$ is $k$-graded in the sense of \cite{Cap-Schichl}. 

Let $G$ be a Lie group with Lie algebra $\g$, and define subgroups $G_{0} \subset P \subset G$ by letting $G_{0}$ (resp. $P$) comprise those elements of $\Ad(G)$ preserving the $\integer$-grading (resp. $\integer$-filtration) of $\g$ (alternatively, $\Ad(G_{0})$ (resp. $\Ad(P)$) is the group of $\integer$-graded (resp. $\integer$-filtered) Lie automorphisms of $\g$). It can be shown that the Lie algebra of $G_{0}$ (resp. $P$) is $\g_{0}$ (resp. $\p$). By construction $G_{0}$ (resp. $P$) covers (possibly trivially) a subgroup of the group $G_{\grad}^{\li}(\g_{-})$ (resp. $G_{\filt}^{\li}(\g_{-})$) of $\integer$-graded (resp. $\integer$-filtered) Lie automorphisms of $\g_{-}$. In general the containment $\Ad(G_{0}) \subset G_{\grad}^{\li}(\g_{-})$ can be proper.

If $\q \subset \g$ is a second parabolic subalgebra, so that $\p \cap \q$ is also parabolic, then the scaling elements $E_{\p}$ and $E_{\q}$ determine on $\g$ a compatible $\ztwo$-grading for which the subspace $\g_{i_{1}, i_{2}}$ is the intersection of the $i_{1}$-eigenspace of $\Ad(E_{\p})$ and the $i_{2}$-eigenspace of $\Ad(E_{\q})$. The subspace $\g_{i_{1}, i_{2}}$ is trivial unless either both $i_{1}$ and $i_{2}$ are non-negative or both are non-positive. It will be convenient to write $\g_{0}$ for the subspace $\g_{0, 0}$.

In what follows these constructions will be used mainly in the special case in which $G = Sp(n, \rea)$ and $P$ is one of the subgroups $P_{kl}$ described in Section \ref{flatmodel}. In particular, viewing the subgroup $P_{12}$ as the intersection $P_{1} \cap P_{2}$, the scaling elements $E_{\p_{1}}$ and $E_{\p_{2}}$ determine on $\g = \spl(n, \rea)$ a compatible $\ztwo$-grading with associated compatibe $\integer$-grading determined by the scaling element $E_{\p_{12}} = E_{\p_{1}} + E_{\p_{2}}$. Fix a basis, $\{f_{\infty}, e_{\infty}, e_{1} \dots , e_{2n-4}, e_{0}, f_{0}\}$, in the symplectic vector space, $(\standrep, \Omega)$, so that $\Omega(f_{\infty}, f_{0}) = 1 = \Omega(e_{\infty}, e_{0})$, $\Omega(f_{\infty}, e_{\infty}) = 0 = \Omega(f_{0}, e_{0})$, and $\Omega(e_{i}, e_{j}) = \omega_{ij}$. When $n > 2$, the associated $\integer$-grading is $4$-step.
%REMOVED FOR SPACE
%In the case $n > 2$, the $\ztwo$-grading $\g = \bigoplus_{K}\g_{K}$ is represented schematically by
%\[\stackrel{\begin{pmatrix} \g_{0,0} &  \g_{1,0} &  \g_{1, 1} & \dots &  \g_{1,1} &  \g_{1,2} &  \g_{2,2} \\  \g_{-1,0} &  \g_{0,0} &  \g_{0,1} & \dots &  \g_{0,1} &  \g_{0,2} &  \g_{1,2}\\  \g_{-1, -1} &  \g_{0, -1} &  \g_{0,0} & \dots &  \g_{0,0} &  \g_{0,1} &  \g_{1,1}\\ \vdots & \vdots & \vdots & \dots & \vdots & \vdots & \vdots \\   \g_{-1,-1} &  \g_{0,-1} &  \g_{0,0} & \dots &  \g_{0,0} &  \g_{0,1} &  \g_{1,1}\\  \g_{-1, -2} &  \g_{0,-2} &  \g_{0, -1} & \dots &  \g_{0, -1} &  \g_{0,0} &  \g_{1,0}\\  \g_{-2,-2} &  \g_{-1, -2} &  \g_{-1, -1} & \dots &  \g_{-1, -1} &  \g_{-1, 0} &  \g_{0,0}  \end{pmatrix}  }{\,\,\,\,\quad\underbrace{\qquad\qquad\qquad\qquad}_{2n-4}}.\]
%The associated $\integer$-grading is $4$-step. 
The most general element in $\g_{-}$ has the form
\begin{equation}\label{gengminus}
\begin{pmatrix} 0 & 0 & 0 & 0 & 0\\
x^{\infty} & 0 & 0 & 0 & 0\\
x^{p} & u^{p} & 0 & 0 & 0\\
x^{0} & u^{0} & -u_{q} & 0 & 0 \\
z & x^{0} & -x_{q} & -x^{\infty} & 0 \end{pmatrix} = \begin{matrix} x^{\infty}t_{-1, 0} + u^{p}a_{p} + u^{0}t_{0, -2} \\ \\+ x^{p}e_{p} + x^{0}t_{-1, -2} + zt_{-2, -2}\end{matrix}
\end{equation}
(where, when convenient, $x_{0} = -x^{\infty}$ and $x_{\infty} = x^{0}$). Here $\g_{0} \simeq \spl(n-2, \rea) \oplus \rea \oplus \rea$ and
\begin{align*}
&\g_{0,-1} = \text{span}\left\{a_{i}\right\},&
&\g_{-1, 0} = \text{span}\left\{t_{-1, 0} \right\},&
&\g_{-1, -1}=\text{span}\left\{e_{i}\right\},\\
&\g_{0, -2}= \text{span}\left\{  t_{0, -2}\right\} ,&
&\g_{-1, -2}=\text{span}\{t_{-1, -2}\},&
&\g_{-2, -2}=\text{span}\{t_{-2, -2} \}.
\end{align*}
%DON'T NEED DATA OF ROOTS
%\begin{align}
%&\g_{0,-1}& &=& &\oplus_{j = 1}^{n-2}\g_{\pm\lambda_{j+2} - \lambda_{2}}& &=& &\text{span}\left\{a_{i}, \, 1\leq i \leq 2n-4 \right \},\\
%&\g_{-1, 0}& & =& &\g_{\lambda_{2} - \lambda_{1}}& &=& &\text{span}\left\{t_{-1, 0} \right\},\\
%&\g_{-1, -1}& &=& &\oplus_{j = 1}^{n-2}\g_{\pm\lambda_{j+2} - \lambda_{1}}& &=& &\text{span}\left\{e_{i}, \, 1 \leq i \leq 2n-4 \right\},\\
%&\g_{0, -2}& &=&  &\g_{-2\lambda_{2}}& &=& &\text{span}\left\{  t_{0, -2}\right\} ,\\
%&\g_{-1, -2}& &=& &\g_{-\lambda_{2} - \lambda_{1}}& & =&  &\text{span}\{t_{-1, -2}\},\\
%&\g_{-2, -2}& &=& &\g_{-2\lambda_{1}}& &=& &\text{span}\{t_{-2, -2} \}.
%\end{align}
As $\g_{-}$ is generated as a Lie algebra by $\g_{-1, 0}$ and $\g_{0, -1}$, the following bracket relations completely determine $\g_{-}$:
\begin{align}
\label{gminusbrackets}
&[a_{i}, a_{j}] = -2\omega_{ij}t_{0, -2},&
&[a_{i}, t_{-1, 0}] = e_{i},&\\
\notag  &[t_{-1, 0}, t_{0, -2}] = -t_{-1, -2},&
&[a_{i}, t_{0, -2}] = 0,&\\
\notag &[a_{i}, e_{j}] = -\omega_{ij}t_{-1, -2},&
&[a_{i}, t_{-1, -2}] = 0,&\\
\notag &[t_{-1,0}, t_{-1, -2}] = -2t_{-2, -2},&
&[t_{-1, 0}, e_{i}] = 0.& 
%&[e_{i}, e_{j}] = -2\omega_{ij}t_{-2, -2}.
\end{align}
(The bracket $[e_{i}, e_{j}] = -2\omega_{ij}t_{-2, -2}$ follows from \eqref{gminusbrackets} by the Jacobi identity).

Direct computation shows that $G_{0} \subset Sp(n, \Rea)$ is isomorphic to $Sp(n-2, \rea) \times GL(1, \rea) \times GL(1, \rea)$ and that the adjoint action of the element $(C_{i}\,^{j}, c, d) \in G_{0}$ on $\g_{-}$ is given by
\begin{align}
\label{adg0action}&\bar{t}_{-2, -2} = c^{-2} t_{-2, -2},& &\bar{t}_{-1, -2} = c^{-1}d^{-1} t_{-1, -2},& &\bar{t}_{0, -2} = d^{-2}t_{0, -2},&\\
\notag & \bar{e}_{i} = c^{-1}C_{i}\,^{j}e_{j},& &\bar{t}_{-1, 0} = c^{-1}dt_{-1, 0},& &\bar{a}_{i} = d^{-1}C_{i}\,^{j}a_{j}.
\end{align}
It follows that the kernel of the projection $G_{0} \to \Ad(G_{0})$ is $(\pm\delta_{i}\,^{j}, \pm 1, \pm 1)$. 

\begin{lemma}\label{z2gisos}
Every element of the group, $G_{\grad}^{\li}(\g_{-})$, of $\integer$-graded Lie automorphism of $\g_{-}$ is in fact $\ztwo$-graded and $G_{\grad}^{\li}(\g_{-})$ is isomorphic to $CSp(n-2, \rea) \times GL(1, \rea)$. $Ad(G_{0})$ is isomorphic to the subgroup, $CSp(n-2, \rea) \times GL^{+}(1, \rea) \subset G_{\grad}^{\li}(\g_{-})$, comprising elements the restriction of which to $\g_{-2, -2}$ is orientation preserving.
\end{lemma}

\begin{proof}
Let $\bar{x}$ denote the image of $x \in \g_{-}$ under a $\integer$-graded linear automorphism. Then there are constants $a$, $b$, $c$, $A_{0}^{i}$, $A_{i}\,^{j}$, $A_{i}\,^{0}$, $d$, $B_{0}\,^{i}$, $B_{i}\,^{j}$, and $B_{i}\,^{0}$ so that $c, d \neq 0$ and the matrices $\begin{pmatrix} a & A_{i}\,^{0} \\ A_{0}\,^{i} & A_{i}\,^{j}\end{pmatrix}$ and $\begin{pmatrix} b & B_{i}\,^{0} \\ B_{0}\,^{i} & B_{i}\,^{j}\end{pmatrix}$ are invertible and
\begin{align*}
&\bar{t}_{-1, 0} = at_{-1, 0} + A_{0}\,^{i}a_{i},& &\bar{a}_{i} = A_{i}\,^{j}a_{j} + A_{i}\,^{0}t_{-1, 0},&& \bar{e}_{i} = B_{i}\,^{j}e_{j} + B_{i}\,^{0}t_{0, -2},& \\
 &\bar{t}_{0, -2} = bt_{0, -2} + B_{0}\,^{i}e_{i},&&\bar{t}_{-1, -2} = c t_{-1, -2},& &\bar{t}_{-2, -2} = d t_{-2, -2} .
\end{align*}
Demanding that the given $\integer$-graded linear automorphism be a Lie automorphism implies that the bracket relations \eqref{gminusbrackets} hold with unbarred elements replaced by barred elements. This yields
\begin{align*}
&d = a^{2}b,& &c = ab,& &A_{i}\,^{0} = A_{0}\,^{i} = 0 = B_{i}\,^{0} = B_{0}\,^{i} ,&\\
&A_{i}\,^{p}A_{j}\,^{q}\omega_{pq} = b\omega_{ij},& &B_{i}\,^{j} = aA_{i}\,^{j},
\end{align*}
so that a $\integer$-graded Lie autormorphism has the form
\begin{align}
\label{zgradaction} &\bar{t}_{-1, 0} = at_{-1, 0},& &\bar{a}_{i} = A_{i}\,^{j}a_{j},& &\bar{e}_{i} = aA_{i}\,^{j}e_{j},& \\
\notag &\bar{t}_{-2, -2} = a^{2}bt_{-2, -2},& &\bar{t}_{-1, -2} = abt_{-1, -2},& &\bar{t}_{0, -2} = bt_{0, -2},&
\end{align}
where $a, b \neq 0$ and $A_{i}\,^{p}A_{j}\,^{q}\omega_{pg} = b\omega_{ij}$. Such a transformation evidently preserves also the $\ztwo$-grading of $\g_{-}$. The relation $A_{i}\,^{p}A_{j}\,^{q}\omega_{pg} = b\omega_{ij}$ shows that $A_{i}\,^{j}$ determines $b$, and evidently the pair $(A_{i}\,^{j}, a) \in CSp(n-2, \rea) \times GL(1, \rea)$. Equations \eqref{zgradaction} show that if the restriction to $\g_{-2, -2}$ of an element of $G_{\grad}^{\li}(\g_{-})$ is orientation preserving, then $b > 0$, in which case the corresponding element of $Ad(G_{0})$ of the form \eqref{adg0action} is that element $(C_{i}\,^{j}, c, d)$ determined by setting $d = b^{-1/2}$, $c = da$, and $C_{i}\,^{j} = caA_{i}\,^{j}$.
\end{proof}
By Theorem 5.3 of \cite{Yamaguchi}, when $n > 2$ it is the case that $\g_{0}$ is the Lie algebra of Lie derivations of $\g_{-}$, and this could be used also to prove Lemma \ref{z2gisos}. %Here a direct proof is given.

Because $\g_{-}$ is nilpotent, the exponential map on $\g_{-}$ is a diffeomorphism onto its image, so that \eqref{gengminus} can be taken as coordinatizing $G/P_{12}$. By the Darboux Theorem for contact manifolds, these coordinates $x^{I}$ and $z$ can be used on a general contact manifold, $(M,H)$, so that $H$ is the kernel of the one-form $\theta = dz + \Omega_{PQ}x^{P}dx^{Q} = dz + x^{\infty}dx^{0} - x^{0}dx^{\infty} + \omega_{pq}x^{p}dx^{q}$. The left-invariant basis for $H_{x,z}$ is $X_{I} = \frac{\partial}{\partial x^{I}} + \Omega_{IP}x^{P}\frac{\partial}{\partial z}$, and the one-forms $dx^{I}$ form a dual basis for $H^{\ast}$. Define coordinates, $a^{I}$, on the fiber $H_{x, z}$ by writing, for $X \in H_{x,z}$, $X = a^{I}X_{I} = dx^{I}(X)$. Suppose $a^{\infty} \neq 0$ and set $u^{\alpha} = \frac{a^{\alpha}}{a^{\infty}}$. Then $x^{I}, z, u^{\alpha}$, are a system of coordinates on the open subset of $\ph$ determined by $a^{\infty} \neq 0$. In local coordinates the left-invariant vector fields, $X_{h}(g) = \frac{d}{dt}_{t = 0} g\cdot \exp{(th)}$, are 
\begin{align*}
&A_{i}  = \tfrac{\partial}{\partial u^{i}} + \omega_{ip}u^{p}\tfrac{\partial}{\partial u^{0}},&
&T_{-1, 0}  = X_{\infty} + u^{\alpha}X_{\alpha},&
&E_{i}  = X_{i} + \omega_{ip}u^{p}X_{0},\\
&T_{0, -2}  = \tfrac{\partial}{\partial u^{0}},&
&T_{-1, -2} = X_{0},&
&T_{-2, -2} = \tfrac{\partial}{\partial z}.&
\end{align*}
For example $X_{a_{i}} = A_{i}$. These vector fields satisfy the bracket relations obtained from \eqref{gminusbrackets} by replacing a lower case letter by a capital letter. The following left-invariant coframe on $\ph$ will be used as a reference coframe:
\begin{align*}
&\theta^{-1, 0} = dx^{\infty},&
&\theta^{i} = du^{i},&
&\theta^{-1, -2} = dx^{0} - u^{0}dx^{\infty} + \omega_{pq}u^{p}dx^{q}\\
&\theta^{-2, -2} = \theta,& 
&\eta^{i} = dx^{i} - u^{i}dx^{\infty}.& 
&\theta^{0, -2} = du^{0} + \omega_{pq}u^{p}du^{q}.& 
\end{align*}
The differentials of these one-forms satisfy algebraic relations dual to the bracket relations \eqref{gminusbrackets} and neatly summarized by the statement that $d\Theta + \Theta \wedge \Theta = 0$ where $\Theta$ is the $\g$-valued one-form defined by:
\begin{align*}
\Theta = \begin{pmatrix} 
 0 & 0 & 0 & 0 & 0\\
\theta^{-1, 0} & 0 & 0 & 0 & 0\\
\eta^{i} & \theta^{i} & 0 & 0 & 0\\
\theta^{-1, -2} & \theta^{0, -2} & -\theta_{j}& 0 & 0\\
\theta^{-2, -2} & \theta^{-1, -2} & -\eta_{j} & -\theta^{-1, 0} &0\\
\end{pmatrix}
\end{align*}
 The one-form $\Theta$ is the local coordinate expression of the Cartan connection for the flat model contact path geometry.

For subbundles $E, F \subset TN$, let $\partial(E, F)_{x}$ denote the subspace of $T_{x}N$ spanned by all the elements of the form $X_{x} + Y_{x} + [Z, U]_{x}$ where $X, Z \in \Gamma(E)$ and $Y, U \in \Gamma(F)$. The pair $(E, F)$ is {\bf regular} if $\dim \partial(E, F)_{x}$ is constant. If $(E, F)$ is regular then $\partial(E, F)$ is a subbundle of $TN$ of constant. $\partial(E, E)$, will be denoted $\partial E$. The {\bf $p$th weak derived system of $E$} is defined inductively by $\partial^{(0)} = E$ and $\partial^{(p)} = \partial(E, \partial^{(p-1)}E)$.

\begin{lemma}\label{bracketgeneratingcontact}
For the canonical filtration, \ref{canonicalcontactfiltration}, of $T\proj(H)$, there hold the following bracket relations:
\begin{align*}
&\partial E = \Hup,& & \partial^{(2)}E = T\proj(H) & &\partial U = V,& \\
& \partial(E, V) = \Hup,&
&\partial(E, U) = E^{\perp},&
&\partial(E^{\perp}, V) = \Hup.&
\end{align*}
\end{lemma}

\begin{proof}
The Darboux theorem implies that any contact manifold is locally equivalent to the flat model contact manifold and, consequently, that $\proj(H)$ equipped with the canonical filtration \ref{canonicalcontactfiltration} is locally equivalent to the flat model. In particular, the subbundles $U$, $V$, $E$, $E^{\perp}$, and $\Hup$, of $T\proj(H)$, are locally always as in the flat case. The left invariant vector fields in the flat model are therefore always a convenient frame in $T\proj(H)$, (an observation which greatly simplifies calculations). The claimed bracket relations now follow immediately from the brackets of $\g_{-}$ described in \eqref{gminusbrackets}.
\end{proof}

\begin{lemma}\label{liftlemma}
Every fiber preserving diffeomorphism $\Phi:\proj(H) \to \proj(H)$ preserving $E$ is the lift of a contactomorphism of $(M, H)$. 
\end{lemma}

\begin{proof}
Because $\Phi$ is fiber preserving, it covers a diffeomorphism, $\phi$, on $M$, and from $\phi \circ \pi = \pi \circ \Phi$ there follows $\Phi^{\ast}(\pi^{\ast}(\theta)) = \pi^{\ast}(\phi^{\ast}(\theta))$ for any contact one-form, $\theta$, on $M$. Since $\Phi$ preserves $E$, $\Phi^{\ast}(\pi^{\ast}(\theta))$ lies in $\ann(E) \subset \ann(V)$. If $A \in \Gamma(V)$ and $B \in \Gamma(E)$ then 
\begin{align*}
d\Phi^{\ast}(\pi^{\ast}(\theta))(A, B) = \phi^{\ast}(d\theta)(\pi_{\ast}(A), \pi_{\ast}(B)) = 0,
\end{align*}
which implies that $\Phi^{\ast}(\pi^{\ast}(\theta))([A, B]) = 0$. Because $\Hup = \partial(V, E)$, this shows $\Phi^{\ast}(\pi^{\ast}(\theta)) \in \ann(\Hup)$, so that $\phi^{\ast}(\theta) \in \ann(H)$ and $\phi$ is a contactomorphism. Composing $\Phi$ with the map induced on $\proj(H)$ by $\phi^{-1}$ shows that to prove the lemma it suffices to consider $\Phi$ that covers the identity map on $(M, H)$. That $\Phi$ preserves $E$ means that for $A \in E_{L}$ there holds $\Phi_{\ast}(L)(A) \in E_{\Phi(L)}$, and so $\pi_{\ast}(L)(\Phi_{\ast}(L)(A)) \in \Phi(L) \subset H_{\pi \circ \Phi(L)} = H_{\pi(L)}$. On the other hand, $\pi_{\ast}(L)(\Phi_{\ast}(L)(A)) = \phi_{\ast}(L)(A) \in L \subset H_{L}$. This forces $\Phi(L) = L$.
\end{proof}

The foliation of $\proj(H)$ is generated locally by the integral curves of a vector field, $X \in E$, transverse to the vertical, $V$, 
\begin{equation}\label{mostgeneral}
X = CT_{-1, 0} + f^{0}T_{0, -2} + f^{i}A_{i}.
\end{equation}
where $C \neq 0$ and $f^{\alpha}$, are functions of $x^{I}, z, u^{I}$. 
%OMITTED FOR SPACE
%Explicitly,
%\begin{equation*}
%X = C\frac{\partial}{\partial x^{\infty}} + Cu^{\alpha}\frac{\partial}{\partial x^{\alpha}} + C(x^{0} - x^{\infty}u^{0} + \omega_{pq}u^{p}x^{q} )\frac{\partial}{\partial z}  + f^{p}\frac{\partial}{\partial u^{p}} + (f^{0} + \omega_{pq}f^{p}u^{q})\frac{\partial}{\partial u^{0}}.
%\end{equation*}
%%%The integral curves satisfy
%%%\begin{align*}
%%%&\dot{x}^{\infty} = C,&  &\dot{x}^{\alpha} = Cu^{\alpha},& &\dot{z} = C(x^{0} - x^{\infty}u^{0} + \omega_{pq}u^{p}x^{q}),&\\
%%%&\dot{u}^{p} = f^{p},& &\dot{u}^{0} = f^{0} + \omega_{pq}f^{p}u^{q}.&
%%%\end{align*}
After some simplification, the integral curves of $X$ satisfy:
\begin{align*}
 &\dot{x}^{\infty} = C, & &\dot{x}^{\alpha} = u^{\alpha}\dot{x}^{\infty},& 
&\dot{z} = - \Omega_{IJ}x^{I}\dot{x}^{J},&
&\dot{u}^{p} = f^{p},& &\dot{u}^{0} = f^{0} + \tfrac{1}{C}\omega_{pq}f^{p}\dot{x}^{q}
\end{align*}
These are the most general equations defining locally a contact path geometry. The specialization $\dot{x}^{\infty} = 1$ determines the parameter $t$ up to translation. Since a translation of $t$ leaves the form of the equations unchanged, it may be assumed that $x^{\infty} = t$. With this specialization the equations admit the obvious simplification to the form \eqref{contactpathodes}. In the case that $n = 2$, so that $M$ is three-dimensional, then equations \eqref{contactpathodes} become
\begin{align}\label{cpath3odes}
&\ddot{x} = f(t, x(t), \dot{x}(t), z(t)),& &\dot{z} = x - t\dot{x}.&
\end{align}

\begin{remark}\label{pdqremark}
The contactomorphism $p = (x^{\infty})^{2} + (x^{0})^{2}$ and $-x^{\infty}\tan{q} = x^{0}$ sends $\theta = dz + \Omega_{PQ}x^{P}dx^{Q}$ to $dz - pdq + \omega_{ij}x^{i}dx^{j}$. The image of the span of $X_{\infty}$ and $X_{0}$ is spanned by $\frac{\partial}{\partial q} + p\frac{\partial}{\partial z}$ and $2\frac{\partial}{\partial p}$. In a region of $\proj(H)$ where $q \neq 0$, let $u^{0} = a^{0}/q$ and $u^{i} = a^{i}/q$. Define vector fields on $\proj(H)$,
\begin{align*}
&E_{i} = X_{i} + \omega_{ip}u^{p}\tfrac{\partial}{\partial p},& &T_{-1, -2} = \tfrac{\partial}{\partial p},& &T_{-2, -2} = \tfrac{1}{2}\tfrac{\partial}{\partial z},& \\
&T_{-1, 0} = \tfrac{\partial}{\partial q} + p\tfrac{\partial}{\partial z} + u^{0}\tfrac{\partial}{\partial p} + u^{p}X_{p},&
\end{align*}
where $X_{i}$, $A_{i}$, and $T_{0, -2}$ are as before. These vector fields satisfy the bracket relations \eqref{gminusbrackets}. The most general contact path geometry is equivalent to the graphs of the integral curves of the vector field $X = CT_{-1, 0} + f^{0}T_{0, -2} + f^{p}A_{p}$, which satisfy the equations
%OMITTEd FOR SPACE
%\begin{align*}
%&\dot{q} = C,& &\dot{p} = Cu^{0},& &\dot{x}^{i} = Cu^{i},&\\
%&\dot{z} = C(p + u^{i}x^{j}\omega_{ij}),& &\dot{u}^{i} = f^{i},& &\dot{u}^{0} = f^{0} + f^{i}u^{j}\omega_{ij}.&
%\end{align*}
The specializations $C = 1$ and $q = t$ give
\begin{align}
&\label{pqspec}\ddot{x}^{i} = f^{i},& &\ddot{p} = f^{0} + f^{u}\dot{x}^{v}\omega_{uv},& &\dot{z} = p + \dot{x}^{u}x^{v}\omega_{uv}.& 
\end{align}
If $n = 2$, so that $(M, H)$ is a three-dimensional contact manifold, \eqref{pqspec} becomes $\dddot{z} = f^{0}(t,z(t),\dot{z}(t), \ddot{z}(t))$. This is the most general third order ordinary differential equation studied by Chern, \cite{Chern}, (though some effort is required to explicitly relate his results and those described here). %See also the discussion of related problems given by E. Cartan in \cite{Cartan-Spanish}.
\end{remark}

\subsection{Canonical Filtration}\label{filteredsection}
A \textbf{$\integer$-filtered vector space} is a vector space, $V$, and a collection of subspaces $\{V^{i} \subset V: i \in \integer\}$ such that $V^{i} \cap V^{j} = V^{\max\{i, j\}}$, $\cap_{i \in \integer} V_{i} = \{0\}$, and $\cup_{i \in \integer}V^{i} = V$. The associated $\integer$-graded vector space is $\graded V = \oplus_{i} V^{i}/V^{i+1}$.  If $V$, $W$ are $\integer$-filtered vector spaces, a linear map $\phi:V \to W$ is $\integer$-filtered if $\phi(V^{i}) \subset W^{i}$ for all $i \in \integer$. A {\bf $\integer$-filtered manifold} modeled on the $\integer$-filtered vector space, $V$, is a smooth manifold, $N$, equipped with a collection of smooth subbundles $T^{i} \subset TN$ such that $\dim T^{i} = \dim V^{i}$, $T^{i} \cap T^{j} = T^{\max\{i, j\}}$, $\cap_{i \in \integer}T^{i} = \{0\}$, and $\cup_{i \in \integer}T^{i} = TN$. $\graded TN = \oplus_{i}U_{i}$ has the structure of a $\integer$-graded vector bundle where $U_{i} = T^{i}/T^{i+1}$. A $\integer$-filtered manifold is \textbf{bracket filtered} if it satisfies $[T^{i}, T^{j}] \subset T^{i + j}$ for all $i, j< 0$.

\begin{lemma}\label{cpathbifiltration}
When $2n -1 > 3$, the choice of a splitting $E = V \oplus W$ determines subbundles $T^{I} \subset T\ph$ defined by
\begin{align*}
&T^{0, -1} = U,& 
&T^{-1, 0} = W,& 
&T^{-1, -1} = \partial(U, W),& \\
&T^{0, -2} = V,& 
&T^{-1, -2} = \Hup,& 
&T^{-2, -2} = T\proj(H),& 
\end{align*}
which satisfy $\rank T^{I} = \dim \f_{I}$, $T^{I} \cap T^{J} = T^{\lub\{I, J\}}$, $\cap_{I \in \ztwo}T^{I} = \{0\}$, and $\cup_{I \in \ztwo}T^{I} = T\ph$. Moroeover, the Lie brackets satisfy the following containments:
\begin{align*}
&[T^{-1,0}, T^{-1, 0}] \subset T^{-1, 0},& &[T^{0, -1}, T^{0, -1}] \subset T^{0, -2},&\\
 &[T^{0, -1}, T^{0, -2}] \subset T^{0, -2},&
&[T^{-1, 0}, T^{0,-1}] \subset T^{-1, -1},&\\
&[T^{-1, 0}, T^{0, -2}] \subset T^{-1, -2},&
&[T^{-1, 0}, T^{-1, -1}] \subset T^{-1, -2},&\\
&[T^{0, -1}, T^{-1, -1}] \subset T^{-1, -2},& &[T^{0, -1}, T^{-1, -2}] \subset T^{-1, -2}.&
\end{align*}
\end{lemma}

\begin{proof}
The content of the first part of the lemma is that $T^{-1, -1} = \partial(U, W)$ has rank $4n-7$, is contained in $T^{-1, -2}$, and intersects $T^{0, -2}$ in $T^{0, -1}$. $W$ is spanned by a vector field of the form $X = T_{-1, 0} + fT_{0, -2} + f^{p}A_{p}$ for some smooth functions $f$, $f^{p}$. There follows 
\begin{equation}\label{xaebracket}
[A_{i}, X] = E_{i} + (A_{i}(f) + 2f_{i})T_{0, -2} + A_{i}(f^{p})A_{p}, 
\end{equation}
and so the $4n-7$ linearly independent vector fields, $A_{i}$, $X$, and $[A_{i}, X]$ span $T^{-1, -1}$. Note that $T_{0, -2}$ is not contained in $T^{-1, -1}$, so the intersection $T^{-1, -1}\cap T^{0, -2}$ is the span of the $A_{i}$, which is $T^{0, -1}$. As $\Hup$ is spanned by $A_{i}$, $T_{-1, 0}$, $T_{0,-2}$, and $E_{i}$, \eqref{xaebracket} shows also $T^{-1, -1} \subset T^{-1, -2}$.

$\partial(U, W) \subset \partial(V, E) = \Hup = T^{-1, -2}$, which implies $[T^{-1, 0}, T^{0, -2}] \subset T^{-1, -2}$. Note that it has been shown that $T^{-1, -1} = \text{span}\{A_{i},  X, E_{i} + (A_{i}(f) + 2f_{i})T_{0, -2}\}$, and that $T^{-1, -1} + T^{0, -2} = \text{span}\{A_{i}, T_{0, -2}, T_{-1, 0}, E_{i}\} = E^{\perp}$. Consequently, to determine $[T^{-1, -1}, T^{-1, 0}]$ it suffices to compute
\begin{equation}\label{axxbracket}
[[A_{i}, X], X] = (df(A_{i}) + 3f_{i})T_{-1, -2} \mod T^{-1, -1} + T^{0, -2},
\end{equation}
which shows $[T^{-1, 0}, T^{-1, -1}] \subset T^{-1, -2}$. Likewise, to verify $[T^{0, -1}, T^{-1, -1}] \subset T^{-1, -2}$ it suffices to compute 
\begin{equation}\label{aaxbracket}
[A_{i}, [A_{j}, X]] = -\omega_{ij}T_{-1, -2} \mod T^{0, -2}. 
\end{equation}
This shows $[T^{0, -1}, T^{-1, -1}] \subset T^{-1, -2}$. The remaining claims are obvious.
\end{proof}
Lemma \ref{cpathbifiltration} shows that $[T^{I}, T^{J}] \subset \sum_{K \nless I + J}T^{K}$, where $K \nless I + J$ means that either $K \geq I + J$ or $K$ is not comparable to $I + J$. However \eqref{axxbracket} of the proof of Lemma \ref{cpathbifiltration} shows that in general $[T^{I}, T^{J}]$ is not a subbundle of $T^{I + J}$ as in general $[T^{-1, 0}, T^{-1, -1}]$ is not contained in $T^{-2, -1} = T^{-1, -1}$.

\begin{corollary}\label{cpathsemiregular}
Associate to $E = V \oplus W$ the $\integer$-filtration 
\begin{align*}
& T^{-1} = U\oplus W,& &T^{-2} = E^{\perp},& &T^{-3} = \Hup,& &T^{-4} = T\proj(H).&
\end{align*}
$T^{-1}$ generates this filtration in the sense that $T^{-i} = \partial(T^{-1}, T^{-i+1})$, and consequently this $\integer$-filtration is a bracket filtration.
\end{corollary}

\begin{proof}
By definition $\partial T^{-1} = \partial(T^{-1}, T^{-1})$ is spanned by linear combinations of $A_{i}$, $X$, $[A_{i}, X]$, and $[A_{i}, A_{j}]$. These span the subspace spanned by $A_{i}$, $T_{-1, 0}$, $T_{0, -2}$, and $E_{i}$, which is $E^{\perp} = T^{-2}$, as was shown in the proof of Lemma \ref{cpathbifiltration}. So $\partial T^{-1} = T^{-2}$. The preceeding shows $\partial(T^{-1}, T^{-2})$ contains $T^{-3} = T^{-1, -2} = [T^{0, -2}, T^{-1, 0}]$. The containments of Lemma \ref{cpathbifiltration}%\eqref{sr3} and \eqref{sr8}
 show that $\partial(T^{-1}, T^{-2}) \subset T^{-3}$, so equality must hold. To show $\partial(T^{-1}, T^{-3}) = T^{-4}$, it suffices to check that $T_{-2, -2}$ lies in $\partial(T^{-1}, T^{-3})$. A standard inductive argument now shows that $[T^{i}, T^{j}] \subset T^{i+j}$ for $i, j < -1$. 
\end{proof}

The constituent subbundles of $\graded TN = \oplus_{i = -1}^{-4}U_{i}$ are:
\begin{align*}
&U_{-1} = U\oplus W,& &U_{-2} = E^{\perp}/(U\oplus W),&
 %= (\partial(U, W)/(U \oplus W)) \oplus (E/(U\oplus W)),&\\
 &U_{-3} = \Hup/E^{\perp},& &U_{-4} = T\proj(H)/\Hup. 
\end{align*}
Setting $U_{I} = T^{I}/(\sum_{J > I}T^{J})$ gives subbundles, $U_{I}$, for which $U_{i} = \oplus_{|I| = i}U_{I}$. For instance $U_{-1, -1} = \partial(U, W)/(U \oplus W)$. The realization of each of the $U_{I}$ as a subbundle of the $U_{i}$ for which $i = |I|$ is obvious, except for $U_{0, -2}$. The point is that $T^{0,-2}$ does not contain $T^{-1}$. However, $T^{0, -2} \subset T^{-2}$, and its image under the projection $T^{-2} \to T^{-2}/T^{-1}$ equals the image of $E$ under the same projection. Consequently the subbundle of $U_{-2}$ corresponding to $U_{0, -2} = V/U$ is $E/(U\oplus W) = (V\oplus W)/(U\oplus W) \simeq V/U$.

\eqref{axxbracket} shows that for any choices of $X$ and $X^{\prime} = CX$, each spanning $W$, and any section $A$ of $U$, the iterated bracket, $[X, [X^{\prime}, A]]$, is a section of $E^{\perp}$, and that its image in the quotient bundle, $\Hup/E^{\perp}$, is well-defined up to scale. Precisely, for any smooth functions, $a, b$ and $g^{p}$, because $X$ and $[X, A]$ are contained in $E^{\perp}$ there holds $[aX, [bX, g^{p}A_{p}]] = abg^{p}[X, [X, A_{p}]] \,\mod E^{\perp}$. Consequently the map $(X, X^{\prime}, A) \to [X, [X^{\prime}, A]] + E^{\perp}$ determines a section, $\ct$, of $\hup/E^{\perp} \tensor S^{2}(W^{\ast})\tensor U^{\ast}$. The Lie bracket induces a map $U_{-1, 0}\times U_{0, -1} = W \times U \to \partial(U, W)/(U\oplus W)= U_{-1, -1}$ which induces an isomorphism $U_{-1, 0} \tensor U_{0, -1} \simeq U_{-1, -1}$, given by $(cX, a^{p}A_{p}) \to ca^{p}[X, A_{p}] + E^{\perp}$. Using this isomorphism the map $U_{-1, 0} \times U_{-1, -1} \to U_{-1, -2}$ induced by the Lie bracket of vector fields can be identified with the section $\ct$.

\begin{definition}\label{vancontor}
The section $\ct$ of $\hup/E^{\perp} \tensor S^{2}(W^{\ast})\tensor U^{\ast}$ is the \textbf{contact torsion} of the contact path geometry. A contact path geometry has \textbf{vanishing contact torsion} if $\ct = 0$.
\end{definition}
Evidently a contact path geometry has vanishing contact torsion if for any choice of $X$ spanning $W$ and any section $A$ of $U$, the iterated bracket, $[X, [X, A]]$ is contained in $E^{\perp}$. Since for three-dimensional contact path geometries, $E^{\perp} = E$ and $U$ is trivial, three-dimensional contact path geometries necessarily have vanishing contact torsion. 

\begin{remark}\label{torsionidentityremark}
For the particular choices of $A_{i}$ and $X$ given above, explicit computation utilizing the fact that $E^{\perp}$ is spanned by $A_{i}$, $T_{0, -2}$, $X$, and $E_{i}$ shows that $\ct = 0$ occurs if and only if $3f_{p} + A_{p}(f^{0}) = 0$. When $\ct = 0$, the contact paths are locally equivalent to the solutions of a system of equations \eqref{contactpathodes} which has now the form \eqref{cpathodenotorsion}.
\end{remark}

\subsection{Canonical Cartan Connections}
The basic terminology about Cartan connections used here follows that of \cite{Cap-Survey} which can be consulted for further references. The Lie algebra homology $H_{\ast}(\pplus; \g)$ is the homology of the complex $(C_{k}(\pplus;\g) = \Lambda^{k}(\pplus)\tensor \g, \partial)$ where 
\begin{align*}
&\partial((p_{1} \wedge \dots \wedge p_{k})\tensor a) = \sum_{i = 1}^{k-1}(-1)^{i}(p_{1}\wedge \dots \wedge  \hat{p}_{i}\wedge \dots \wedge p_{k})\tensor [p_{i}, a] +\\& \sum_{i < j}(-1)^{i+j}([p_{i}, p_{j}]\wedge p_{1} \wedge \dots \wedge \hat{p}_{i} \wedge \dots \wedge \hat{p}_{j}\wedge \dots \wedge p_{k})\tensor a.
\end{align*}
The curvature function, $\kappa$, of a $(\g, P)$ Cartan connection on $\pi:\G \to N$, is the $P$-equivariant $\Lambda^{2}((\g/\p)^{\ast})\tensor \g$-valued function on $\G$ defined by $\kappa(u)(h_{1}, h_{2}) = \K_{u}(\eta_{u}(h_{1}, \eta_{u}(h_{2}))$. The Killing form on $\g$ gives an isomorphism $(\g/\p)^{\ast} \simeq \pplus$, and so $\kappa(u)$ may be viewed as an element of $C_{2}(\pplus; \g)$. The group $P$ acts naturally on the complex $C_{\ast}(\pplus; \g)$, and $\partial$ is $P$-equivariant, and so $P$ acts also on the homologies $H_{\ast}(\pplus;\g)$. The scaling element $E_{\p}$ determined by $P$ determines a filtration (grading) on any $P$-module ($G_{0}$-module). A $(\g, P)$ Cartan connection is called \textbf{regular} if there vanishes for all $p \leq 0$ the homogeneity $p$ part, $\kappa^{(p)}$, of the its curvature function, where $\kappa^{(p)}$ is the part of $\kappa$ lying in $(\f_{i}^{\ast} \wedge \f_{j}^{\ast})\tensor \f_{i+j+p}$. The action of $\pplus$ on $H_{\ast}(\pplus; \g)$ is trivial, and so $H_{\ast}(\pplus; \g)$ can be viewed as a $G_{0}$-module equipped with the $\integer$-grading $H_{k}(\pplus; \g) = \oplus_{l}H^{(l)}_{k}(\pplus; \g)$ induced by the scaling element $E_{\p} \in \g_{0}$. Consult Section 3 of \cite{Calderbank-Diemer} for details. Hence there can be considered the associated bundles $\G\times_{P}C_{k}(\pplus; \g)$ and $\G\times_{P}H_{k}(\pplus;\g)$, and $\kappa$ may be regarded as a section of $\G\times_{P}C_{k}(\pplus; \g)$. If $\partial(\kappa(u)) = 0$ for every $u \in \G$, then $\eta$ is called \textbf{normal}; in this case $\kappa$ determines a section of $\G\times_{P}H_{k}(\pplus; \g)$.

Call \textbf{compatible} a $(\g, P_{12})$ Cartan connection, $\eta$, on $\pi:\adaptedframe \to N$ if it induces the canonical filtration \eqref{canonicalcontactfiltration} in the sense that 
\begin{align*}
&\pi_{\ast}(\eta^{-1}(\f_{0, -1})) = U,&
&\pi_{\ast}(\eta^{-1}(\f_{0, -2} + \f_{-1, 0})) = E,&
&\pi_{\ast}(\eta^{-1}(\f_{-1, -2})) = \Hup,&\\ &\pi_{\ast}(\eta^{-1}(\f_{0, -2})) = V,& &\pi_{\ast}(\eta^{-1}(\f_{0, -2} + \f_{-1, -1})) = E^{\perp}. 
\end{align*}
A compatible $\eta$ induces a contact path geometry with $W = \pi_{\ast}(\eta^{-1}(\f_{-1, 0}))$. 

\begin{theorem}\label{cpathnormal}
A contact path geometry is induced by a compatible, regular $(\g, P_{12})$ Cartan connection on $\G$ if and only if it has vanishing contact torsion. Among the compatible $(\g, P_{12})$ Cartan connections on $\pi:\adaptedframe \to \proj(H)$ inducing a contact torsion free contact path geometry there is a unique isomorphism class of regular, normal Cartan connections.
\end{theorem}

Let $\g$ be a semisimple Lie algebra and $\p \subset \g$ a parabolic subalgebra as in the discussion at the beginning of Section \ref{localdescriptioncontactpath}. If a smooth manifold, $N$, has an $\integer$-bracket filtration, the graded vector bundle, $\graded TN = \bigoplus_{i < 0}(T^{i}/T^{i+1})$, acquires from the Lie bracket of vector fields on $N$ a fiberwise algebraic Lie bracket defined as follows. Given $X \in T^{i}_{x}/T^{i+1}_{x}$ and $Y \in T^{j}_{x}/T^{j+1}_{x}$, choose vector fields $\bar{X}\in \Gamma(T^{i})$ and $\bar{Y} \in \Gamma(T^{j})$ such that the images in $T^{i}_{x}/T^{i+1}_{x}$ and $T^{j}_{x}/T^{j+1}_{x}$, respectively, of $\bar{X}_{x}$ and $\bar{Y}_{x}$, respectively, equal $X$ and $Y$, respectively. The image, $[X, Y]_{x}$, in $T^{i+j}/T^{i+j+1}$ of $[\bar{X}, \bar{Y}]$ does not depend on the choice of $\bar{X}$ and $\bar{Y}$, and so defines the desired algebraic Lie bracket. This bracket makes $\graded TN$ a bundle of $\integer$-graded Lie algebras. The Lie algebra $\graded T_{x}N$ usually is called the \textbf{symbol algebra at $x$}.

If the $\integer$-filtration of $TN$ is modeled on $\g_{-}$, the fiber over $x \in N$ of the graded frame bundle, $\frameN$, of $\graded TN$, comprises all graded linear isomorphisms $u:\g_{-} \to \graded T_{x}N$. This is a principal bundle with structure group $GL_{\grad}(\g_{-})$, the group of $\integer$-graded linear automorphisms of $\g_{-}$. The adjoint representation identifies $\Ad(G_{0})$ with a subgroup of the group $GL_{\grad}(\g_{-})$. Because the Lie bracket on $\g_{-}$ is $\Ad(G_{0})$-invariant, a reduction of the structure group of $\frameN$ to a subgroup, $\Ad(G_{0})$, induces a fiberwise Lie bracket on $\graded TN$ defined for $X, Y \in \graded T_{x}N$ by $\{X, Y\}_{x} = u([u^{-1}(X), u^{-1}(Y)])$, where $u$ is any element of the reduced bundle lying over $x$. A reduction of the structure group of $\gr TN$ to $\Ad(G_{0})$ is called \textbf{regular} if $\{\,,\,\}_{x} = [\,,\,]_{x}$ for all $x \in N$. Precisely, regard  the Lie bracket on $\g_{-}$ as a distinguished element $c_{\g} \in \Lambda^{2}(\g_{-}^{\ast})\tensor \g_{-}$. The fiberwise Lie bracket $[\,,\,]$ on $\graded TN$ induces a function $c:\frameN \to \Lambda^{2}(\g_{-}^{\ast})\tensor \g_{-}$ defined by $c(u)(a, b) = u^{-1}([u(a), u(b)]_{\pi(u)})$ for $a, b \in \g_{-}$. A reduction of $\frameN$ to $\Ad(G_{0})$ is regular if the the restriction of $c$ to the reduced bundle is constant and equal to $c_{\g}$.

For simple $G$, a theorem equivalent to Theorem \ref{prolongationtheorem} is due to Tanaka, \cite{Tanaka-equivalence}. The version stated here, covering the case that $G$ is semisimple, is  Corollary 3.23 in \v{C}ap- Schichl's \cite{Cap-Schichl}.
\begin{theorem}[A. \v{C}ap - H. Schichl, T. Morimoto, N. Tanaka]\label{prolongationtheorem}
Let $G$ be a finite dimensional semisimple Lie group with $\integer$-graded Lie algebra, $\g$. Let $N$ be a $\integer$-filtered manifold modeled on $\g_{-}$. If $H_{1}^{(p)}(\pplus; \g) = \{0\}$ for all $p > 0$, the isomorphism classes of regular reductions to the structure group $\Ad(G_{0})$ of the graded frame bundle are in bijection with isomorphism classes of principal $P$-bundles over $N$ endowed with a regular, normal Cartan connection. 
\end{theorem}
\noindent
With respect to the cohomological hypothesis of Theorem \ref{prolongationtheorem}, recall 
\begin{theorem}[Yamaguchi, Theorem 5.3 of \cite{Yamaguchi}]\label{cohomologicaltheorem}
Let $\g$ be a simple graded Lie algebra over $\rea$ such that $\g_{p} = [\g_{p+1}, \g_{-1}]$ for $p < -1$. Then $H_{1}^{(p)}(\pplus; \g) = 0$ for all $p > 0$ except when $\g$ is isomorphic with $(A_{n}, \{\alpha_{1}\})$, $(C_{n}, \{\alpha_{1}\})$ or the split real form of either of these.
\end{theorem}
\noindent
Here the notation $(A_{n}, \{\alpha_{1}\})$ (resp. $(C_{n}, \{\alpha_{1}\})$) indicates the grading of $sl(l+1, \rea)$ (resp. $\spnr$) determined by marking the first node of the Dynkin diagram $A_{n}$ (resp. $C_{n})$). The two exceptional series correspond to projective and contact projective structures. On the other hand, the series $(C_{l}, \{\alpha_{1}, \alpha_{2}\})$ correspond to the contact path geometries, and Yamaguchi's theorem shows that these do satisfy the hypothesis of Theorem \ref{prolongationtheorem}. As is explained below in Section \ref{curvaturesection} the relevant homologies are straightforward to compute using Kostant's version of the Bott-Borel-Weil theorem.

It had been hoped that Theorem \ref{prolongationtheorem}, or one of its analogues (see \cite{Cap-Schichl}, \cite{Doubrov-Komrakov-Morimoto}, \cite{Morimoto}, \cite{Morimoto-survey}, \cite{Tanaka-equivalence})
%%OMITTED FOR SPACE
% (The most general such theorem is the one stated by Morimoto in \cite{Morimoto-survey}) could be applied to the general contact path geometry. 
That there should be a difficulty is apparent already with the solution of the equivalence problem for the contact projective structures found in \cite{Fox-cproj}. A contact projective structure with contact torsion is induced by a $(\g, P_{1})$ Cartan connection which is regular but not normal. By Proposition 2.4 and Theorem 2.7 in \cite{Cap}, the $(\g, P_{12})$ Cartan connection inducing the contact path geometry determined by such a contact projective structure is also not normal (and will in general not be regular). Other natural examples of non-regular Cartan connections arise via \v{C}ap's correspondence space construction, see e.g. Section 4.6 of \cite{Cap}.

Tanaka and Morimoto formulated differently the hypothesis of a regular $\Ad(G_{0})$ reduction of the structure group of $\graded TN$. Section 4.4 of \cite{Cap-Schichl} explains carefully the correspondence between the different formulations. For a statement of a theorem like Theorem \ref{prolongationtheorem} see the survey \cite{Doubrov-Komrakov-Morimoto} (in fact the statement there is more general). See also the survey \cite{Yamaguchi-Yatsui}. Theorem 2.3 of Morimoto's \cite{Morimoto} is stronger than Theorem \ref{prolongationtheorem}. However its hypotheses are still too strong to apply to a contact path geometry with torsion. On the other hand, Morimoto's general theory does apply in principle to the situation of contact path geometries with torsion, although it does not associate to them anything so nice as a Cartan connection.

\begin{remark}\label{structuredataremark}
In practice the data of a regular reduction of the graded frame bundle may be constructed explicitly by exhibiting a local frame in $\graded TN$ in which the Lie brackets agree with the Lie brackets in $\g_{-}$. More formally, such a local frame is a local section, $u:N \to \frameN$, of the graded frame bundle such that the function $c(u)$ on $N$ is constant and equal to $c_{\g_{-}}$. If $\tilde{u}$ is any other local section of $\frameN$ such that $c(\tilde{u}) = c_{\g_{-}}$, then over the intersection in $N$ of the domains of $u$ and $\tilde{u}$ there exists a function, $b$, taking values in $Gl_{\grad}^{\li}(\g_{-})$ such that $\tilde{u} = u \circ b$. It follows that the $Gl_{\grad}^{\li}(\g_{-})$ orbit of any such frame gives a reduction of the structure group to $Gl_{\grad}^{\li}(\g_{-})$. (This is closely related to the discussion following Proposition 3.5.2 in \cite{Morimoto}). Precisely, the reduced bundle comprises exactly those $u \in \frameN$ such that $c(u) = c_{\g}$. 

In those cases in which the Lie algebra of $Gl_{\grad}^{\li}(\g_{-})$ is $\g_{0}$, the existence of a frame in which the Lie brackets agree with those in $\g_{-}$ essentially (i.e. modulo some discrete topological data) already determines a reduction of the structure group to $\Ad(G_{0})$. In the setting of contact path geometries (with vanishing contact torsion) by virtue of Lemma \ref{z2gisos} such a procedure works. In general, the data of a filtration regular of type $\g_{-}$ determines only a reduction to $G_{\grad}^{\li}(\g_{-})$ of the structure group of $\graded TN$, and some extra data is needed to determine the required $\Ad(G_{0})$ reduction. An example of this sort is given by the case $SL(n+1, \rea)/P_{1 k+1}$ corresponding to the $k$-path geometries modeled on the manifolds of flags of points and $k$-planes in projective space. For finite-dimensional real simple Lie algebras admitting a compatible $\integer$-grading it is the case that $\gl_{\grad}^{\li}(\g_{-}) = \g_{0}$ if and only if $H_{1}^{(p)}(\pplus; \g) = 0$ for all $p > 0$. The cases for which this fails are enumerated in Theorem 5.3 of \cite{Yamaguchi}. 

On the other hand, to prove that there exists no regular $\Ad(G_{0})$ reduction of the graded frame bundle, it suffices to show that there is no choice of local frame in $\graded TN$ for which the Lie brackets agree with those of $\g_{-}$.
\end{remark}

\begin{proof}[Proof of Theorem \ref{cpathnormal}]
Suppose given a contact path geometry induced by a compatible, regular $(\g, P_{12})$ Cartan connection, $\eta$, on $\pi:\G \to \ph$, and let $\kappa$ be the curvature function of $\eta$. 
%Compatibility implies $W = \pi_{\ast}(\eta^{-1}(\f_{-1, 0}))$, $U = \pi_{\ast}(\eta^{-1}(\f_{0, -1})$, and $V = \pi_{\ast}(\eta^{-1}(\f_{0, -2}))$. Using these 
Using the compatibility of $\eta$ it is straightforward to check that regularity implies that $\partial(U, W) = \pi_{\ast}(\eta^{-1}(\f_{-1, -1}))$, which can then be used with regularity to show that $\ct = 0$. 

If there exists a regular $\Ad(G_{0})$ reduction of the graded frame bundle, $\frameN$, of $\graded TN$, then for any section $u:N \to \redframe$ of the $\Ad(G_{0})$ reduced bundle, there must hold $[u(a_{i}), u(t_{-1, 0})] = u(e_{i})$ and $[u(e_{i}), u(t_{-1, 0})] = 0$. Together these imply $[[u(a_{i}), u(t_{-1, 0})], u(t_{-1, 0})] = 0$. Consequently, for any $a \in \g_{0,-1}$ and any $b \in \g_{-1, 0}$, there must hold $[[u(a), u(b)], u(b)] = 0$. The most general section of $U_{0, -1} = U$, has the form $g^{p}A_{p}$, and so the most general brackets of a section of $U_{0, -1}$ with sections of $U_{-1, 0} = W$ have the form
\begin{align*}
& [g^{p}A_{p}, CX] = Cg^{p}(E_{p} + (df(A_{p}) +2f_{p})T_{0, -2}) \mod T^{-1, 0}+ T^{0, -1},\\
& [[g^{p}A_{p}, CX], CX] = C^{2}g^{p}(df(A_{p}) + 3f_{p})T_{-1, -2} \mod T^{-1, -1} + T^{0, -2},
\end{align*}
and from these it is evident that if $3f_{i} \neq -df(A_{i})$ then there is no choice of a frame in $\graded TN$ for which the Lie brackets agree with those of $\g_{-}$.

On the other hand, if $\ct = 0$ then there holds $3f_{i} = - df(A_{i})$. Define a local section, $u$, of $\frameN \to \ph$ by $u(a_{i}) = A_{i}$, $u(t_{0, -2}) = T_{0, -2} + T^{-1}$, $u(t_{-1, -2}) = T_{-1, -2} + T^{-2}$, $u(t_{-2, -2}) = T_{-2, -2} + T^{-3}$, and
\begin{align*}
&u(t_{-1, 0}) = T_{-1, 0} + f^{0}T_{0, -2} + f^{p}A_{p},
&u(e_{i}) = E_{i} - f_{i}T_{0, -2} + T^{-1}.& 
\end{align*}
Direct computation confirms that $c(u) = c_{\g_{-}}$, i.e. that the Lie brackets in this frame agree with those of $\g_{-}$. As explained in Remark \ref{structuredataremark} this $u$ determines a $Gl_{\grad}^{\li}(\g_{-})$ reduction of $\frameN$ which evidently depends only on the contact structure and the splitting $E = V \oplus W$. By Lemma \ref{z2gisos}, the subgroup of $Gl_{\grad}^{\li}(\g_{-})$ inducing orientation-preserving automorphisms of $\g_{-2, -2}$ is isomorphic to $\Ad(G_{0})$. A co-orientation of the underlying contact structure is an orientation of $U_{-2, -2}$, and a reduction of the structure group of $\frameN$ to $\Ad(G_{0})$ is effected by considering only graded frames, $u:\g_{-} \to \graded T_{\pi(u)}N$ mapping a fixed orientation on $\g_{-2, -2}$ to the given orientation on $U_{-2, -2}$. Alternatively this canonical $\Ad(G_{0})$ reduction of $\frameN$ is determined by considering the $\ztwo$-graded frames in $\graded TN$ consistent with the orientation of $U_{-2, -2}$. This reduction is regular and Theorem \ref{cohomologicaltheorem} shows that the cohomological hypotheses of Theorem \ref{prolongationtheorem} are satisfied. By Theorem \ref{prolongationtheorem} there exists a principal $P_{12}$-bundle and on it a unique isomorphism class of regular, normal Cartan connection inducing the given contact path geometry. By definition any $(\g, P_{12})$-Cartan connection determines an absolute parallelism on the supporting $P_{12}$-bundle, and so this bundle is topologically trivial. Hence, given any concrete realization of a trivial $P_{12}$-bundle, this isomorphism class of Cartan connections may be viewed as living on said bundle. As a connection representing this isomorphism class must induce the filtration \eqref{canonicalcontactfiltration}, there is a Cartan connection on $\adaptedframe$ representing the given isomorphism class.
\end{proof}

\subsubsection{Split Quaternionic Contact Structures}\label{splitquatstructure}
Isotropic Grassmannian structures are parabolic geometries modeled on $Q_{k}$. This section treats the cases $2 \leq k \leq n-1$, with special emphasis on the case $k = 2$ (the simpler case $k = n$ corresponding to the Lagrangian Grassmannian is omitted).

\begin{definition}\label{isograssdefined}
For $n \geq k +1$, call a codimension $k(k+1)/2$ bracket generating distribution, $T^{-1}N$, on a $(2k(n-k) + k(k+1)/2)$-dimensional smooth manifold, $N$, an \textbf{isotropic Grassmannian structure} if there is a regular reduction to $\Ad(G_{0})$ of the structure group of the associated graded vector bundle $\gr TN$. The structures arising in the case $k = 2$ are called \textbf{split quaternionic contact structures}. 
\end{definition}

By Theorem \ref{prolongationtheorem}, an isotropic Grassmannian structure induces a canonical regular normal $(\g, P_{k})$ Cartan connection on a $P_{k}$ principal bundle over $N$. Definition \ref{isograssdefined} is unsatisfactory from a geometric point of view because it is not stated in terms of readily identifiable (e.g. tensorial) geometric structures on $N$. Proposition \ref{aisoprop} shows that such a reduction has a straighforward geometric meaning. Proposition \ref{splitredprop} gives a further reformulation of the special case of split quaternionic contact structures.

\begin{definition}
An \textbf{almost isotropic Grassmannian structure} on $N$ consists of a codimension $\binom{k+1}{2}$ bracket generating distribution, $C$; a rank $2(n-k)$ symplectic vector bundle, $(V, \Omega)$; and a rank $k$ vector bundle $L$ such that $C \simeq V \tensor L^{\ast}$ and $TN/C \simeq S^{2}(L^{\ast})$. 
\end{definition}

\begin{proposition}\label{aisoprop}
The data of an almost isotropic Grassmannian structure determines a regular $\Ad(G_{0})$ reduction of the structure group of $\gr TN$ (and so determines an isotropic Grassmannian structure) if and only if the $S^{2}(L^{\ast})$-valued symplectic structure on $C$ induced by the identification $C \simeq V \tensor L^{\ast}$ agrees with the $TN/C$-valued $2$-form $\bomega$ determined by the Lie bracket of vector fields on $C$ (and defined in section \ref{flatmodel}).
\end{proposition}

\begin{proof}

The parabolic $P_{k}$ of section \ref{flatmodel} determines a $\integer$-grading of $\g = \spnr$ such that the graded Lie subalgebra $\g_{-} = \g_{-1}\oplus \g_{-2}$ is isomorphic to that spanned by the $X_{i}\,^{\al}$ and the $V^{\al\be}$. For $k \geq 2$, by Theorem 5.3 of \cite{Yamaguchi} (or Theorem 5.1 of \cite{Yamaguchi-Yatsui}) $\g$ is the algebraic prolongation of $\g_{-}$, which implies that $\g_{0}$ comprises graded Lie derivations of $\g_{-}$. This means the group $Gl_{\grad}^{\li}(\g_{-})$ has Lie algebra $\g_{0}$. The subalgebra $\g_{0}$ is isomorphic to the product $\smp(n-k, \rea) \times \gl(k, \rea)$. With this information in hand it is not hard to check that $Gl_{\grad}^{\li}(\g_{-})$ is isomorphic to the quotient of the product $CSp(n-k, \rea) \times GL(k, \rea)$ by the $\rea^{\times}$ action $r\cdot (C, D) = (rC, r^{-1}D)$, where the action on $\g_{-}$ of $(C_{i}\,^{j}, D^{\al}\,_{\be}) \in CSp(n-k, \rea) \times GL(k, \rea)$ (with $C_{i}\,^{p}C_{j}\,^{q}\omega_{pq} = c\omega_{ij}$, $c^{n-k} = \det C_{i}\,^{j}$) is given by $(X_{i}\,^{\al}, V^{\al\be})\to (C_{i}\,^{p}D^{\al}\,_{\si}X_{p}\,^{\si}, cD^{\al}\,_{\si}D^{\be}\,_{\ga}V^{\si\ga})$. (Checking this statement directly %without using Theorem 5.3 of \cite{Yamaguchi} 
is tedious). The Lie algebra of $\Ad(G_{0})$ is also $\g_{0}$, and it can be checked easily that $\Ad(G_{0})$ is the quotient, $Sp(n-k, \rea)\cdot GL(k, \rea)$, of $Sp(n-k, \rea) \times GL(k, \rea)$ by the equivalence $(C, D) \sim  (-C, -D)$.
%(or, what is the same, the quotient of $CSp_{0}(n-k, \rea) \times GL(k, \rea)$ by the equivalence $(C, D) \sim (rC, r^{-1}D)$ for $r \in \rea^{\times}$, where the subscript indicates the connected component of the identity).

If $N$ has an almost isotropic Grassmannian structure, the vector bundle $\gr TN = (V \tensor L^{\ast})\oplus S^{2}(L^{\ast})$ equipped with the fiberwise algebraic Lie bracket $[(u\tensor a, r), (v \tensor b, s)] = (0, -\Omega(u, v)(a \tensor b + b \tensor a))$ is fiberwise isormorphic as a Lie algebra to $\g_{-}$, so its structure group admits a reduction to $Gl_{\grad}^{\li}(\g_{-})$. Note that by the remarks in the preceeding paragraph, any linear automorphism of $\gr TN$ is induced by a conformal symplectic automorphism of $(V, \Omega)$ and a linear automorphism of $L$. To effect a further reduction of the structure group to $\Ad(G_{0})$ it suffices to impose some condition forcing the conformal symplectic automorphism of $(V, \Omega)$ to be symplectic.

A $\Lambda^{2}(L^{\ast})$-valued symmetric two-form, $\bg$, on $V \tensor L^{\ast}$ is defined by $\bg(u\tensor a, v \tensor b) = \Omega(u, v)a \wedge b$. View $\bg$ as a section of $T = S^{2}((V\tensor L^{\ast})^{\ast})\tensor \Lambda^{2}(L^{\ast})$. The action of $CSp(n-k, \rea) \times GL(k, \rea)$ on $T$ does not descend to an action of $Gl_{\grad}^{\li}(\g_{-})$ on $T$, but it does give an action of $Gl_{\grad}^{\li}(\g_{-})$ on the oriented projectivization, $\proj^{+}(T)$. Consequently, a reduction of the structure group of $\gr TN$ is effected by restricting to the subgroup of $Gl_{\grad}^{\li}(\g_{-})$ the action of which on $\proj^{+}(T)$ fixes the image %in $\proj^{+}(T)$ 
of $\bg$. Explicit calculation shows that this amounts to considering only elements of $Gl_{\grad}^{\li}(\g_{-})$ induced by elements of $CSp_{0}(n-k, \rea) \times Gl(k, \rea)$, i.e. this effects the desired reduction to $\Ad(G_{0})$.

Finally, the statement that the $\Ad(G_{0})$ reduction is regular if and only if the $S^{2}(L^{\ast})$-valued symplectic structure on $C$ induced by the identification $C \simeq V \tensor L^{\ast}$ agrees with the $TN/C$-valued $2$-form $\bomega$ %determined by the Lie bracket of vector fields on $C$ 
is just a rephrasing of the definition of regularity.
\end{proof}

In the particular case of $k = 2$, split quaternionic contact structures may be formulated in terms of purely tensorial data analogous to the data used to define quaternionic contact structures in \cite{Biquard}.

\begin{definition}
A codimension $3$ bracket generating distribution, $C$, on the $(4n-5)$-dimensional ($n \geq 3$) smooth manifold, $\paths$, is an \textbf{almost split quaternionic contact structure} if there is given on $C$ a split signature conformal structure, $[g]$, and locally on $\paths$ there exist on $C$ endomorphisms $E, F, J$ satisfying the relations $E^{2} = F^{2} = I = -J^{2}$ and $E \circ F = J$, such that for any representative $g \in [g]$, the tensors $g(J\cdot, \cdot)$, $g(E\cdot, \cdot)$ and $g(F\cdot, \cdot)$ are skew-symmetric.
\end{definition} 
An almost split quaternionic contact structure on $\paths$ determines a reduction to $\Ad(G_{0})$ of the structure group of the associated graded bundle $\graded T\paths$. Define on $C$ an $\sltwor$-valued two-form, $\Bar{\Omega}_{\al}\,^{\be}$, as the restriction to $C$ of  
\begin{align}\label{splitsquatform}
&\Bar{\Omega}_{\be}\,^{\al} = \begin{pmatrix}g(F\cdot, \cdot) & g(E\cdot, \cdot) - g(J\cdot, \cdot)\\ g(E\cdot,\cdot) + g(J\cdot,\cdot) & -g(F\cdot, \cdot)  \end{pmatrix},
\end{align}
and set $\Psi = \det \Bar{\Omega}$ on $C$. Working locally so that $C$ and $TN/C$ are trivial, and using a section of $\Lambda^{2}(TN/C)$ to lower the index of $\Bar{\Omega}_{\be}\,^{\al}$, there can be defined from $\Bar{\Omega}$ on $\gr TN$ an algebraic Lie bracket that is fiberwise isomorphic to $\g_{-}$. This is enough to determine a reduction of the structure group to $Gl_{\grad}^{\li}(\g_{-})$, and the conformal structure $[g]$ can be used as in the proof of Proposition \ref{aisoprop} to determine a further reduction to $\Ad(G_{0})$.

\begin{proposition}\label{splitredprop}
The reduction to $\Ad(G_{0})$ of the associated graded bundle, $\graded T\paths$, of an almost split quaternionic contact structure on $\paths$ is regular (and so is a split quaternionic contact structure) if and only if the restriction to $C$ of the $T\paths/C$-valued $2$-form $\bomega$ is determined as in \eqref{splitsquatform} for any choice of $\eta \in \Lambda^{2}(L^{\ast})$ and any choice of trivializations.  
\end{proposition}

It will be shown in Section \ref{contacttorsion} that sometimes the space of contact paths of a contact path geometry admits a split quaternionic contact structure. 

\subsubsection{Existence of $(\g, P_{12})$ Cartan connections inducing a given contact path geometry}

As will be sketched now, any contact path geometry is induced by some $(\g, P_{12})$-Cartan connection. The requisite $P_{12}$ principal bundle $\G$ is built directly from the filtration \eqref{canonicalcontactfiltration}, and the existence on it of Cartan connections inducing the given contact path geometry follows straightforwardly from the existence of a canonical regular, normal Cartan connection associated to any path geometry subordinate to the given contact path geometry.

As explained in \cite{Fox-cproj}, a co-oriented contact structure and a choice of a square-root of the bundle of positive contact one-forms determines on $M$ a rank $2n$ symplectic vector bundle, the tractor bundle $\tractor$, and the pullback of $\tractor$ over $N$ is a symplectic vector bundle with a $2$-step filtration determined (see below) in part by the containment $V \subset E$. The projectivization of $\tractor$ depends only on the given co-oriented contact structure; $\G$ is defined to be the bundle of filtered symplectic frames in the pulled back tractor bundle. Any path geometry subordinate to the given contact path geometry determines a $(\sll(2n, \rea), \q_{12})$-Cartan connection on a $Q_{12}$-principal bundle ($SL(2n, \rea)/Q_{12}$ is the manifold of $(1, 2)$-flags) over $\ptm$, and the restriction to $\ph$ of said bundle admits a reduction to $\G$. The restriction to $\G$ of the Cartan connection determined by the subordinate path geometry decomposes as the sum of a $(\g, \p_{12})$-Cartan connection on $\G$ and an $\m$-valued one-form on $\G$, where $\m$ is the Killing orthogonal complement of $\spl(n, \rea)$ in $\sll(2n, \rea)$ (essentially this follows from the observation that the contact projective tractor bundle may be indentified with the projective tractor bundle equipped with a symplectic structure). This shows that there exist always Cartan connections on $\G$ inducing the given filtration of $TN$. The difficult matter is to identify normalizations on the curvatures of such Cartan connections that distinguish among them a unique representative (or a useful family of representatives).

For completeness, the construction of $\adaptedframe$ is described, though many details are omitted. Choose a square-root, $\rho:\form \to M$, of the bundle of positive contact one-forms on the co-oriented contact manifold, $(M, H)$. A tautological one-form, $\alpha$, on $\form$, is defined by $\alpha_{s}(X) = s^{2}(\rho_{\ast}(X))$. Let $\delta_{r}:\form \to \form$ denote the principal $\rea^{\times}$-action of $r$. Denote by $\emf[k]$ the real line bundle associated to $\form$ by the representation of $\rea^{\times}$ on $\rea$ given by $s \cdot t = s^{-k}t$. Sections $h \in \Gamma(\emf[k])$ are in canonical bijection with functions $\tilde{h}:\form \to \rea$ of homogeneity $k$, in the sense that $\tilde{h}(\delta_{r}(p)) = r^{k}\tilde{h}(p)$. 
%%OMITTED FOR SPACE
%It will be convenient to write $\emf^{i}[k] = \emf[k]\tensor H$. 
The Jacobi identity shows that the Lie algebra, $\vect(\form)$, of vector fields on $\form$ is graded as a Lie algebra by the homogeneity degree, $\vect(\form) = \oplus_{k \in \integer}\vect_{k}(\form)$. The \textbf{tractor bundle}, $\tractor$, is the rank 2n quotient of $T\form$ by an $\rea^{\times}$ action, $P_{r}$, on $T\form$ covering $\delta_{r}$ and leaving invariant $\vect_{-1}(\form)$, and defined by $P_{r}(Z) = r^{-1}\delta_{r^{-1}}\,^{\ast}(Z)$. By construction, the space of sections, $\Gamma(\tractor)$, is identified with $\vect_{-1}(\form)$. Because the symplectic form, $\Omega = d\alpha$, on $\form$ has homogeneity $2$, it descends to $\tractor$ as a fiberwise symplectic form. The $P_{r}$-invariant filtration $T^{2}\form \subset T^{1}\form \subset T\form$ descends to a filtration $\tractor^{2} \subset \tractor^{1} \subset \tractor$. 
%%OMITTED FOR SPACE
%There is a canonical isomorphism of graded vector bundles $\graded\tractor \simeq \emf[-1]\oplus\emf^{i}[-1]\oplus\emf[1]$.
The tractor bundle depends only on $\form$ and the co-oriented contact structure on $M$. Because $\form^{2n}$ is canonically identified with the bundle of non-vanishing volume forms on $M$, and under this identification $\Omega^{n}$ is identified with a constant multiple of the canonical volume form on $\form^{2n}$, forgetting the symplectic structure and the subbundle $\tractor^{1}$, recovers the usual projective tractor bundle (as described, e.g., in \cite{Cap-Gover}).

Let $\ctractor = \pi^{\ast}(\tractor)$, $\ctractor^{4} = \pi^{\ast}(\tractor^{2})$, and $\ctractor^{1} = \pi^{\ast}(\tractor^{1})$. The symplectic form on $\tractor$ lifts to a symplectic form, also denoted $\Omega$, on $\ctractor$. Recall the filtration $E/V \subset E^{\perp}V \subset \Hup/V \simeq \Hpull$ from Section \ref{contactpathbasics}. There are uniquely determined subbundles $\ctractor^{2}$ and $\ctractor^{3}$ such that $\ctractor^{4} \subset \ctractor^{3} \subset \ctractor^{2} \subset \ctractor^{1}$ and such that, under the isomorphism $\ctractor^{1}/\ctractor^{4} \simeq \Hpull \tensor \pi^{\ast}(\emf[-1])$, $\ctractor^{3}/\ctractor^{4}$ corresponds to $E/V \tensor \pi^{\ast}(\emf[-1])$ and $\ctractor^{2}/\ctractor^{4}$ corresponds to $(E^{\perp}/V)\tensor \pi^{\ast}(\emf[-1])$. Though it is not necessary, $\ctractor^{2}$ and $\ctractor^{3}$ can be defined directly. To each $L \in \proj(H)$, associate a subspace $L^{\prime} \subset \tractor_{\pi(L)}$ as follows. View $L$ as a subspace of $H_{\pi(L)}$ and choose $s \in \form$ lying over $\pi(L)$. Because $P_{r}$ is a rescaling of the differential of the fiber dilations, the image, $L^{\prime}$, in $\tractor_{\pi(L)}$ of the subspace of $T_{s}\form$ projecting onto $L$, does not depend on the choice of $s$. Letting $\ctractor^{3}_{L} = \{t \in \ctractor_{L}: t \in L^{\prime} \subset \tractor_{\pi(L)}\}$ defines $\ctractor^{3} \subset \ctractor^{1}$. $\ctractor^{2}$ is defined to be the skew complement in $\ctractor$ of $\ctractor^{3}$. The bundle of filtered symplectic frames, $\pi:\adaptedframe \to \proj(H)$, in $\ctractor$ is a principal $P_{12}$ bundle the fiber of which over $L \in \proj(H)$ comprises all filtration preserving symplectic linear isomorphisms $u:\standrep \to \ctractor_{L}$. The projectivized bundle $\proj(\ctractor)$ does not depend on the choice of $\form$. The corresponding bundle of filtered projective symplectic frames, $\hat{\pi}:\pframe \to \proj(H)$ is the principal $\bar{P}_{12}$ bundle the fiber of which over $L \in \proj(H)$ comprises all filtration-preserving projective symplectic linear isomorphisms $\proj(\standrep) \to \proj(\ctractor_{L})$. If there had been considered instead the pullback to $\ptm$ of the tractor bundle, $\tractor$, there could have been constructed similarly a principal $Q_{12}$ bundle the fiber of which over $L \in \ptm$ comprises all filtration preserving linear isomorphisms, and from this description it is evident that the restriction to $\ph$ of this bundle admits a canonical reduction to $P_{12}$ which coincides with the bundle $\adaptedframe$ just described.

\section{Curvature of Contact Path Geometries}\label{curvaturesection}
This section has as its goal the geometric interpretation of various torsions and curvatures of a contact path geometry of dimension at least five. Most of the results are proved directly using only the explicit bracket relations from section \ref{prolongsection}. However at a certain point it is necessary to use the full machinery of parabolic geometries. Precisely, when the contact torsion vanishes there is a secondary contact torsion which measures the failure of the generating bundle $W$ to be the characteristic subbundle of $\partial(U, W)$. When this secondary contact torsion vanishes there is determined on the space of contact paths a canonical codimension three multicontact structure. Actually, much more is true - namely this multicontact structure is a split quaternionic contact structure. Proving this claim requires the machinery of harmonic curvature components of a regular, normal parabolic geometry. 

The curvature function, $\kappa$, of a $(\g, P)$ Cartan connection, $\eta$, on $\G$ is viewed as a section of the associated bundle with fiber $\Lambda^{2}(\g_{-}^{\ast})\tensor \g$. When $\eta$ is normal, there is a harmonic Hodge theory allowing computation of the cohomology, $H^{\ast}(\g_{-}; \g) \simeq H_{\ast}(\pplus; \g)$, via Kostant's version of the Bott-Borel-Weil theorem, and it makes sense to speak of the harmonic part of the curvature, $\kappah$. For details the reader should consult \cite{Cap-Survey}, \cite{Cap}, or \cite{Calderbank-Diemer}. Tanaka proved that for regular, normal Cartan connections $\kappah$ is the complete obstruction to local flatness in the sense that $\kappa  = 0$ if and only if $\kappah = 0$. In Theorem 5.10 \cite{Calderbank-Diemer} D. Calderbank and T. Diemer improved this substantially by constructing a natural linear differential operator, $L$, such that $L(\kappah) = \kappa$. Using this A. \v{C}ap gave the following characterization of correspondence spaces.
\begin{theorem}[Theorem 3.3 and Proposition 3.3 of \cite{Cap-Survey}]\label{harmoniccorrespondence}
Let $Q \subset P$ and let $\eta$ be a regular, normal $(\g, Q)$ Cartan connection on the $Q$ principal bundle $\G \to N$. The restriction to a sufficiently small open neighborhood of $N$ of the parabolic geometry $(\eta, \G)$ is isomorphic to the correspondence space of some $(\g, P)$ Cartan connection if and only if $i(x)\kappah = 0$ for all $x \in \p/\q$.
\end{theorem}
Here if $Q \subset P$ a $(\g, P)$ Cartan connection $\eta$ on a $P$ principal bundle $\G \to N$ induces on $\G \to \G/Q$ a $(\g, Q)$ Cartan connection, and the pair $(\eta, \G \to \G/Q)$ is called the correspondence space of the original parabolic geometry. Thereom \ref{harmoniccorrespondence} will be used here to identify when a contact path geometry with vanishing contact torsion is the correspondence space of a split quaternionic contact structure, or, equivalently, when the space of contact paths of a contact path geometry with vanishing contact torsion admits a split quaternionic contact structure.

\subsection{Harmonic Curvature Components}
Representations of a parabolic $P$ will be represented using labeled Dynkin diagrams with crossed nodes; the number over a given node represents the coefficient of the corresponding fundamental weight in the highest weight of the representation. Since $\g = \spnr$ is a split real form the cohomologies $H^{2}(\p_{+}; \g)$ look exactly as they do in the complex case and are computable using the algorithmic implementation of Kostant's version of the Bott-Borel-Weil theorem explained in \cite{Baston-Eastwood}. First one finds the Hasse diagram for the parabolic $P$ and then one uses this to compute a BGG resolution of the representation of $P$ induced by the adjoint representation on $\g$. The harmonic Hodge theory shows $H^{\ast}(\g_{-}; \g) \simeq H_{\ast}(\pplus; \g)$ which, via the usual duality between homology and cohomology, is isomorphic to $H^{\ast}(\pplus, \g^{\ast})^{\ast}$, and $H^{\ast}(\pplus, \g^{\ast})$ can be computed directly using Kostant's theory. Here only $\g = \spnr$ and $\g = \sll(n+1, \rea)$ are used, and in either case the adjoint representation is irreducible and isomorphic to its dual, so one could compute just as well $H^{\ast}(\pplus, \g)$.

In computing the harmonic curvature components it is helpful to know the relation between a representation of $P$ and its dual (this is particular to both $\g$ and the parabolic $P$):
\begin{center}
\parbox[l][1.7\totalheight][c]{158pt}{$\xbbbbdbdual{\mu_{1}}{\lambda_{1}}{\lambda_{2}}{\lambda_{3}}{}{\lambda_{n-2}}{\lambda_{n-1}}$} 
\parbox[r][1.7\totalheight][c]{108pt}{$\xbbbbdb{\hat{\mu}_{1}}{\lambda_{1}}{\lambda_{2}}{\lambda_{3}}{}{\lambda_{n-2}}{\lambda_{n-1}}$} 
\end{center}
where $\hat{\mu}_{1} = - \mu_{1} - 2\sum_{i = 1}^{n-1}\lambda_{i}$;

\begin{center}
\parbox[l][1.7\totalheight][c]{158pt}{$\bxbbbdbdual{\lambda_{1}}{\mu_{2}}{\lambda_{3}}{\lambda_{4}}{}{\lambda_{n-2}}{\lambda_{n-1}}$} 
\parbox[r][1.7\totalheight][c]{108pt}{$\bxbbbdb{\lambda_{1}}{\hat{\mu}_{2}}{\lambda_{3}}{\lambda_{4}}{}{\lambda_{n-2}}{\lambda_{n-1}}$} 
\end{center}
where $\hat{\mu}_{2} = - \lambda_{1} - \mu_{2} - 2\sum_{i = 3}^{n-1}\lambda_{i}$; and when $n \geq 2$,

\begin{center}
\parbox[l][1.7\totalheight][c]{158pt}{$\xxbbbdbdual{\mu_{1}}{\mu_{2}}{\lambda_{1}}{\lambda_{2}}{}{\lambda_{n-3}}{\lambda_{n-2}}$} 
\parbox[r][1.7\totalheight][c]{108pt}{$\xxbbbdb{-\mu_{1}}{\hat{\mu}_{2}}{\lambda_{1}}{\lambda_{2}}{}{\lambda_{n-3}}{\lambda_{n-2}}$} 
\end{center}
\noindent
where $\hat{\mu}_{2} = -\mu_{2} - 2\sum_{i=1}^{n-2}\lambda_{i}$. The last takes the form
\begin{center}
\parbox[l][1.7\totalheight][c]{78pt}{$\xxdual{\mu_{1}}{\mu_{2}}$} 
\parbox[r][1.7\totalheight][c]{108pt}{$\xx{-\mu_{1}}{-\mu_{2}}$} when $n = 2$.
\end{center}

The following tables list the harmonic curvature components for regular normal $(\g, P)$ Cartan connections for $\g = \spnr$ and $P$ which is $P_{1}$, $P_{2}$, or $P_{12}$. With each component are listed its homogeneity and the subspace in which it occurs. In the case $P = P_{12}$ the two scaling elements $E_{\p_{1}}, E_{\p_{2}} \in \g_{0}$ induce a $\ztwo$-grading of any $\g_{0}$-module, in particular of $H_{\ast}(\pplus; \g)$; in this case the unique harmonic representative of the homology class in $H_{2}(\pplus; \g)$ of the curvature function of a regular, normal $(\g, P_{12})$ Cartan connection may be decomposed into $\ztwo$-graded pieces $\kappa^{(K)} \in \Lambda^{2}(\g_{I}^{\ast}\wedge \g_{J}^{\ast})\tensor \g_{I+J+K}$, and the resulting $\ztwo$-homogeneity is recorded in the table. Because $E_{\p_{12}} = E_{\p_{1}} + E_{\p_{2}}$, the ordinary homogeneity is recovered by summing the $\ztwo$-homogeneity. The $n = 2$ case, which has a different character, is omitted. The homogeneities are computed by considering the actions of the scaling elements $E_{\p_{1}}$ and $E_{\p_{2}}$ described in Section \ref{flatmodel}. Precisely, the action of the scaling element $E_{\p_{i}}$ is computed by taking the inner product of the vector comprising the coefficients over the nodes of the Dynkin diagram with the $i$th column of the inverse of the Cartan matrix corresponding to $C_{n}$ (this is explained in \cite{Slovak}).\\

\smallskip
\begin{center}
I. $P = P_{1}$, $n > 2$. 
\end{center}
\begin{center}
\begin{tabular}{|l|c|r|}
\hline
Component of $H_{2}(\pplus; \g)$ & Homogeneity & Contained in\\
\hline
\parbox[c][1.7\totalheight][c]{108pt}{$\xbbbbdb{-1}{2}{1}{0}{}{0}{0}$} & $2$& $\Lambda^{2}(\g_{-1}^{\ast})\tensor \g_{0}$\\
\hline
\end{tabular}
\end{center}
\smallskip

\begin{center}
II. $P = P_{2}$, $n > 3$. 
\end{center}
\begin{center}
\begin{tabular}{|l|c|r|}
\hline
Component of $H_{2}(\pplus; \g)$ & Homogeneity & Contained in\\
\hline
\parbox[c][1.7\totalheight][c]{108pt}{$\bxbbbdb{4}{-3}{0}{1}{}{0}{0}$} & $0$& $\Lambda^{2}(\g_{-1}^{\ast})\tensor \g_{-2}$\\
\hline
\parbox[c][1.7\totalheight][c]{108pt}{$\bxbbbdb{0}{-3}{4}{0}{}{0}{0}$} & $2$& $\Lambda^{2}(\g_{-1}^{\ast})\tensor \g_{0}$\\
\hline
\end{tabular}
\end{center}
\smallskip

\begin{center}
III. $P = P_{2}$, $n = 3$. 
\end{center}
\begin{center}
\begin{tabular}{|l|c|r|}
\hline
Component of $H_{2}(\pplus; \g)$ & Homogeneity & Contained in\\
\hline
\parbox[c][1.7\totalheight][c]{108pt}{$\bxbbbdbt{5}{-3}{1}$} & $1$& $(\g_{-1}^{\ast}\wedge \g_{-2}^{\ast})\tensor \g_{-1}$\\
\hline
\parbox[c][1.7\totalheight][c]{108pt}{$\bxbbbdbt{0}{-3}{4}$} & $2$& $\Lambda^{2}(\g_{-1}^{\ast})\tensor \g_{0}$\\
\hline
\end{tabular}
\end{center}
\smallskip

\begin{center}
IV. $P = P_{12}$, $n > 3$. 
\end{center}
\begin{center}
\begin{tabular}{|l|c|r|}
\hline
Component of $H_{2}(\pplus; \g)$ & Homogeneity & Contained in\\
\hline
\parbox[c][1.7\totalheight][c]{108pt}{$\xxbbbdb{-4}{1}{0}{1}{}{0}{0}$} & $(-2, 0)$& $\Lambda^{2}(\g_{0, -1}^{\ast})\tensor \g_{-2, -2}$\\
\hline
\parbox[c][1.7\totalheight][c]{108pt}{$\xxbbbdb{5}{-4}{1}{0}{}{0}{0}$} & $(2, -1)$& $(\g_{-1, 0}^{\ast}\wedge \g_{-1, -1}^{\ast})\tensor \g_{0, -2}$\\
\hline
\parbox[c][1.7\totalheight][c]{108pt}{$\xxbbbdb{0}{-3}{4}{0}{}{0}{0}$} & $(1, 2)$& $(\g_{0, -1}^{\ast}\wedge \g_{-1, -1}^{\ast})\tensor \g_{0}$\\
\hline
\end{tabular}
\end{center}
\smallskip

\begin{center}
V. $P = P_{12}$, $n = 3$. 
\end{center}
\begin{center}
\begin{tabular}{|l|c|r|}
\hline
Component of $H_{2}(\pplus; \g)$& Homogeneity & Contained in\\
\hline
\parbox[c][1.7\totalheight][c]{108pt}{$\xxbbbdbt{-5}{2}{1}$} & $(-2, 1)$& $(\g_{0, -1}^{\ast}\wedge \g_{0, -2}^{\ast})\tensor \g_{-2, -2}$\\ 
\hline
\parbox[c][1.7\totalheight][c]{108pt}{$\xxbbbdbt{5}{-4}{1}$} & $(2, -1)$& $(\g_{-1, 0}^{\ast}\wedge \g_{-1, -1}^{\ast})\tensor \g_{0, -2}$ \\
\hline
\parbox[c][1.7\totalheight][c]{108pt}{$\xxbbbdbt{0}{-3}{4}$} & $(1, 2)$&  $(\g_{0, -1}^{\ast}\wedge \g_{-1, -1}^{\ast})\tensor \g_{0}$\\
\hline
\end{tabular}
\end{center}
\smallskip

%OMITTED BECAUSE NOT USED
%\begin{center}
%III. $n = 2$ 
%\end{center}
%\smallskip
%\begin{center}
%\begin{tabular}{|l|c|r|}
%\hline
%Component of $H_{2}(\pplus; \g)$& Homogeneity & Contained in\\
%\hline
%\parbox[c][1.7\totalheight][c]{28pt}{$\xx{6}{-3}$} & $(3, 0)$& $(\g_{-1, 0}^{\ast}\wedge \g_{-1, -1}^{\ast})\tensor \g_{1, -1}$\\ 
%\hline
%\parbox[c][1.7\totalheight][c]{28pt}{$\xx{-4}{5}$} & $(1, 3)$& $(\g_{0, -1}^{\ast} \wedge \g_{-2, -1}^{\ast}) \tensor \g_{-1, 1}$\\
%\hline
%\end{tabular}
%\end{center}
%\smallskip

\noindent
When the second node is crossed, the case $n = 3$ looks a little different from the case $n > 3$ because in this case the second node receives a double bond. In either case, the regularity of the Cartan connection implies that any component of homogeneity less than or equal to $0$ vanishes, and so for contact torsion free contact path geometries there are two possibly non-vanishing harmonic curvature components. A harmonic curvature component of homogeneity $K$ must lie in a subspace of $\Lambda^{2}(\g_{-}^{\ast})\tensor \g$ containing the highest weight vector corresponding to this component. The subspace in which lies the component of homogeneity $K$ can be found by calculating the weights of the components $(\g_{I}^{\ast}\wedge \g_{J}^{\ast}) \tensor \g_{I+J+K}$. Often many of the possibilities may be eliminated by simpler considerations. For instance, in cases IV. and V. it is easy to check that a component in $(\g_{I}^{\ast}\wedge \g_{J}^{\ast})\tensor  \g_{I + J + (2, -1)}$ must vanish unless $I = (-1, 0)$ and $J = (-1, -1)$, as all other choices of $I, J < (0, 0)$ give $I + J + (2, -1)$ for which $\g_{I + J +(2, -1)} = \{0\}$. In cases IV. and V. the harmonic curvature component of homogeneity $(1, 2)$ must lie in a subspace of $\Lambda^{2}(\g_{-}^{\ast})\tensor \g$ containing the highest weight vector corresponding to this component. To determine in what subspace lies the component of homogeneity $(1, 2)$ it is necessary to calculate the weights of the components $(\g_{I}^{\ast}\wedge \g_{J}^{\ast}) \tensor \g_{I+J+(1, 2)}$; doing so one finds that only the subspace given by $I = (0, -1)$ and $J = (-1, -1)$ has the required weight. Note that in this case the computations work in the $n = 3$ case just as in the $n > 3$ case. As another example, in case II. the harmonic curvature component of homogeneity $2$ lies in the subspace $(\g_{\lambda_{3}-\lambda_{2}}^{\ast}\wedge \g_{\lambda_{3}-\lambda_{1}}^{\ast})\tensor \g_{2\lambda_{3}}$ where the $\lambda_{i}$ are the standard basis of the dual of a Cartan subalgebra of $\g = \spnr$ which is diagonal in the standard representation.

In making computations it is useful to recall that the highest weight of the representation dual to a given representation is the negative of the lowest weight of the given representation.

\begin{proposition}\label{normalimpliesctf}
The contact path geometry induced on $(M, H)$ by a normal, compatible $(\g, P_{12})$ Cartan connection, $\eta$, on $\pi:\G \to \ph$ has vanishing contact torsion. 
\end{proposition}
\begin{proof}
Tables IV. and V. above show that in all cases there is exactly one possibly non-zero harmonic curvature component of non-positive $\integer$-homogeneity. It is claimed that compatibility forces this component to vanish, so that all possibly non-zero harmonic curvature components have positive $\integer$-homogeneity. Arguing as in the proof of Corollary 4.10 of \cite{Cap-Schichl} by decomposing the Bianchi identity by $\integer$-homogeneities, it is straightforward to check that a normal Cartan connection for which there vanish all harmonic curvature components of non-positive $\integer$-homogeneity must be regular. The proof is completed by observing that a contact path geometry induced by a regular $(\g, P_{12})$ Cartan connection has vanishing contact torsion. To show that compatibility of $\eta$ implies the vanishing of the harmonic curvature component of non-positive $\integer$-homogeneity it is necessary to consider separately the cases $2n -1 > 5$ and $2n-1 = 5$, since in the former case this component has $\ztwo$-homogeneity $A = (-2, 0)$, whereas in the latter case it has $\ztwo$-homogeneity $A = (-2, 1)$. In either case if $x \in \g_{I} \subset \g_{-}$ and $y \in \g_{J} \subset \g_{-}$, then $\kappa^{(A)}(x, y) = \kappa_{A + I + J}(x, y)$ will be automatically $0$ unless $i_{1} + j_{1} - 2 \geq -2$, which, since $i_{1}, j_{1} \leq 0$ implies $i_{1} = 0 = j_{1}$. In the case $2n-1 > 5$, $\kappa_{A + I + J}(x, y)$ can be non-zero if and only if $-2 \leq i_{2} + j_{2} < 0$, which implies $i_{2} = -1 = j_{2}$. For $x, y \in \g_{0, -2}$, $\kappa^{(-2, 0)}(x, y) = \kappa_{-2, -2}(x, y) = \eta_{-2, -2}([\eta^{-1}(x), \eta^{-1}(y)])$. By comptability $\pi_{\ast}(\eta^{-1}(x)) \in \Gamma(U)$ for $x \in \g_{0, -1}$; since $V = \partial(U)$ this implies that $\pi_{\ast}([\eta^{-1}(x), \eta^{-1}(y)]) \in \Gamma(V)$ from which follows $\eta_{-2, -2}([\eta^{-1}(x), \eta^{-1}(y)]) = 0$. The case $2n -1 = 5$ similarly reduces to consideration of $\kappa^{(-2, 1)}(x, y)$ with $x \in \g_{0, -1}$ and $y \in \g_{0, -2}$.
\end{proof}

\smallskip
\subsection{Geometric Meaning of Torsions and Curvatures}\label{contacttorsion}
A contact path geometry is \textbf{admissible} if its space of contact paths, $\paths(M, H, W)$, admits a smooth structure so that the canonical projection $\nu:\ph \to \paths$ is a smooth submersion. Locally every contact path geometry is admissible.

If $\ct$ is not $0$ there holds only $\partial(W, \partial(U, W)) \subset \Hup$. If $2n -1 \geq 5$, by definition $\ct = 0$ if and only if  $\partial(W, \partial(U, W)) \subset E^{\perp}$, and there is in general an obstruction to $\partial(W, \partial(U, W)) \subset \partial(U, W)$. When $\ct = 0$ the Lie bracket of vector fields induces a map, $W \times \partial(U, W)/(U \oplus W) \to E^{\perp}/\partial(U, W)$, which, using the isomorphism $E^{\perp}/\partial(U\oplus W) \simeq V/U$ discussed in section \ref{filteredsection}, can be viewed as a section of $(U_{0, -1}^{\ast}\wedge U_{-1, -1}^{\ast})\tensor U_{0, -2}$, and will be called the \textbf{secondary contact torsion}. The cumbersome terminology `secondary contact torsion' reflects that this invariant is directly analogous to the invariant of a path geometry called the torsion by D. Grossman, \cite{Grossman}. Using the isomorphism $U_{-1, 0} \tensor U_{0, -1} \simeq U_{-1, -1}$ induced by the Lie bracket of vector fields, the secondary contact torsion can be identified with the section of $S^{2}(U_{-1, 0}^{\ast})\tensor U_{0, -1}^{\ast} \tensor U_{-1, -1}  = S^{2}(W^{\ast})\tensor U^{\ast}\tensor E^{\perp}/\partial(U, W)$ defined by $(X, X^{\prime}, A) \in W \times W \times U \to [X, [X^{\prime}, A]] + \partial(U, W) \in E^{\perp}/\partial(U, W) \simeq V/U$.

Call a pair, $\phi^{-2, -2}, \phi^{-1, -2}$, of linearly independent sections of $\ann(E^{\perp})$ a \textbf{filtered frame} in $\ann(E^{\perp})$ if $\phi^{-2, -2}$ spans $\ann(\Hup)$. Any filtered frame in $\ann(E^{\perp})$ has the form
\begin{align*}
&\phi^{-2, -2} = a \theta^{-2, -2},& &\phi^{-1, -2} = b\theta^{-1, -2} + c\theta^{-2, -2},& &a \neq 0 \neq b.
\end{align*}

\begin{lemma}\label{symplecticformlemma}
For any filtered frame in $\ann(E^{\perp})$ the restriction to $E^\perp$ of the associated two-form $d\phi^{-2, -2} + d\phi^{-1, -2}$ gives a symplectic structure on $E^{\perp}$ for which $V$ is a Lagrangian subbundle of $E^{\perp}$ and $E$ is the skew complement in $E^{\perp}$ of $U$.
\end{lemma}

\begin{proof}
The form $\beta = \phi^{-2, -2} + \phi^{-1, -2}$ has the form $r\theta^{-2, -2} + s\theta^{-1, -2}$ with $s \neq 0$. It is claimed that the restriction to $E^{\perp} \times E^{\perp}$ of $d\beta$ is a symplectic structure. The vector fields $A_{i}$, $T_{0, -2}$, $T_{-1, 0}$, and $E_{i}$ span $E^{\perp}$. Direct computation using the flat reference coframe of Section \ref{localdescriptioncontactpath} gives
\begin{align}
\label{intmultintermed} & i(A_{i})d\beta = -s\eta_{i} \mod \ann(E^{\perp}),& \\
\notag &i(T_{0, -2})d\beta = -s\theta^{-1, 0}\mod \ann(E^{\perp}),& \\\notag &i(T_{-1, 0})d\beta = s\theta^{0, -2}\mod \ann(E^{\perp}),& \\
\notag &i(E_{i})d\beta = -2r\eta_{i} - s\theta_{i}\mod \ann(E^{\perp}).
\end{align}
As $-\eta_{i}(E_{j}) = \omega_{ij} = -\theta_{i}(A_{j})$, there holds $d\beta(A_{i}, E_{j}) = s\omega_{ij}$. Likewise, $d\beta(T_{0, -2}, T_{-1, 0}) = -s$. With the non-degeneracy of $\omega_{ij}$ these show that $d\beta$ is non-degenerate when restricted to $E^{\perp}$. Evidently the non-degeneracy holds for any choices of $s \neq 0 $ and $r$. That $V$ is Lagrangian in $E^{\perp}$ and that the skew complement in $E^{\perp}$ of $U$ is $E$ are both evident from \eqref{intmultintermed}.
\end{proof}

Note that if $\tilde{\beta} = f\beta + g\theta^{-2, -2}$ with $f \neq 0$, then $d\tilde{\beta} = fd\beta + gd\theta^{-2, -2} \mod \ann(E^{\perp})$, so the resulting symplectic structure on $E^{\perp}$ does depend on the choice of filtered frame. Note also that if $E = V \oplus W$ is any contact path geometry then it is evident from \eqref{intmultintermed} that $U \oplus W$ is a Lagrangian subbundle of $E^{\perp}$ with respect to the symplectic structure on $E^{\perp}$ determined by any choice of filtered frame in $\ann(E^{\perp})$. If $Q$ is a subbundle of $E^{\perp}$ transverse to $V$ that is Lagrangian with respect to the symplectic structure determined on $E^{\perp}$ by any choice of filtered frame and such that $U \subset Q \subset E$, then a splitting $Q= U \oplus W$ determines a contact path geometry (and every contact path geometry arises in this way). However, such a splitting is not uniquely determined, i.e. there can be $\tilde{W}$ such that $U \oplus \tilde{W} = Q = U \oplus W$.

\begin{proposition}\label{torsionfreeqprop}
Among the contact path geometries determining a given Lagrangian $Q$ such that $U \subset Q \subset E$ and $Q$ is transverse to $V$, there is a unique one which has vanishing contact torsion.
\end{proposition}

\begin{proof}
Let $Q = U \oplus W$ and let $X$ span $W$. Let $\tilde{X} = X + g^{p}A_{p}$ span $\tilde{W}$, so that $U \oplus \tilde{W} = U \oplus W$. Assume without loss of generality that $X = T_{-1, 0} \mod V$. By \eqref{aaxbracket} there are functions $\tau_{i}$ such that $[[A_{i}, X], X] = \tau_{i}T_{-1, -2} \mod E^{\perp}$. Then
\begin{align*}
&[[A_{i}, \tilde{X}], \tilde{X}]= [[A_{i}, X], X] - g^{p}[A_{p}, [A_{i}, X]] + 2g_{i}[T_{0, -2}, X] \mod E^{\perp}\\
& = (\tau_{i} + 3g_{i})T_{-1, -2} \mod E_{\perp}
\end{align*}
so that the contact path geometry determined by $\tilde{W}$ has vanishing contact torsion if and only if $g^{i} = -(1/3)\tau^{i}$. 
\end{proof}

\begin{proposition}\label{twosymp}
The skew complement of $W$ in $E^{\perp}$ with respect to the symplectic structure determined on $E^{\perp}$ by the choice of a filtered frame in $\ann(E^{\perp})$ does not depend on the symplectic structure chosen. This skew complement equals $\partial(U, W)$ if and only if $\ct = 0$.
\end{proposition}

\begin{proof}
Let $\beta =  \phi^{-2, -2} + \phi^{-1, -2} = r\theta^{-2, -2} + s\theta^{-1, -2}$ with $s \neq 0$, and let $X = T_{-1, 0} + f^{0}T_{0, -2} + f^{p}A_{p}$ span $W$. Using \eqref{intmultintermed} one finds $i(X)d\beta = s(\theta^{0, -2} - f^{0}\theta^{-1, 0} - f^{i}\eta_{i}) \mod \ann(E^{\perp})$, and from this it is evident that the kernel of the restriction to $E^{\perp}$ of $i(X)d\beta$ depends neither on the choice of $X$ nor on the choice of filtered frame. Now $d\beta([A_{i}, X], X) = -\beta([[A_{i}, X], X]) = s\theta^{-1, -2}([[A_{i}, X], X])$, and by definition of $\ct$ this vanishes if and only if $\ct = 0$. As $\partial(U, W)$ is spanned by $A_{i}$, $X$, and $[A_{i}, X]$, this shows that the skew complement of $W$ in $E^{\perp}$ equals $\partial(U, W)$ if and only if $\ct = 0$.
\end{proof}

Proposition \ref{twosymp} shows how to determine synthetically from a Lagrangian $U \subset Q \subset E$ the splitting $Q = U \oplus W$ such that $W$ has $\ct = 0$; simply define $W$ to be the skew complement in $E^{\perp}$ of $\partial(U, Q)$.

Recall that for a smooth distribution, $C$, on a smooth manifold, $N$, the characteristic system, $\chsub(C)$ is defined for each $x \in N$ by
\begin{align*}
\chsub_{x}(C) = \{V \in C_{x}: i(V)d\omega \in \ann(C) \,\,\forall \,\, \omega \in \ann(C)\}.
\end{align*}
When $\chsub(C)$ has constant non-zero rank, it is called the \textbf{characteristic subbundle} of $C$. The characteristic system of $C$ measures the degeneracy of $C$. When $\chsub_{x}(C) = \{0\}$ for every $x \in N$, call $C$ \textbf{maximally non-degenerate}.

\begin{proposition}\label{chsubprop}
A contact path geometry with vanishing contact torsion has vanising secondary contact torsion if and only if $W = \chsub(\partial(U, W))$.
\end{proposition}
\begin{proof}
If $\Pi = 0$, $\ann(\partial(U, W))$ is spanned by $\theta^{-2, -2}$, $\theta^{-1, -2}$, and $\mu = i(X)d\theta^{-1, -2}$ for any $X$ spanning $W$. Since by Proposition \ref{twosymp} $\partial(U, W)$ is the skew complement of $W$ in $E^{\perp}$ with respect to the symplectic structure defined on $E^{\perp}$ by $d(r\theta^{-2, -2} + s\theta^{-1, -2})$, the characteristic system at any point of $\ph$ of $\partial(U, W)$ must be at most one-dimensional. As always $\lie_{X}\theta^{-2, -2} = i(X)d\theta^{-2, -2}$ and $\lie_{X}\theta^{-1, -2} = i(X)d\theta^{-1, -2}$ are in $\ann(\partial(U, W))$, it suffices to check that $\lie_{X}\mu = i(X)d\mu$ lies in $\ann(\partial(U, W))$. This is true if and only if $\partial(W,\partial(U, W)) \subset \partial(U, W)$, and the obstruction to the latter is by definition the secondary contact torsion.
\end{proof}

\begin{corollary}[Corollary of Proposition \ref{chsubprop}]
Let $(M, E = V\oplus W)$ be an admissible contact path geometry of dimension at least five with $\ct = 0$ and vanishing secondary contact torsion. Then $\paths$ admits a maximally non-degenerate codimension $3$ bracket generating distribution, $C$, so that the projection $(\ph, \partial(U, W)) \to (\paths, C)$ is a submersion in the category of smooth manifolds equipped with a distribution.
\end{corollary}

Suppose given an admissible contact path geometry induced by a regular, normal $(\g, P_{12})$ Cartan connection, $\eta$, on $\G \to \ph$. Since by the discussion in Section \ref{splitquatstructure} a split quaternionic contact structure on $\paths$ is the same thing as a regular, normal $(\g, P_{2})$ Cartan connection on $\paths$, it makes sense to say that the local space of paths, $\paths$, of an admissible contact path geometry \textbf{admits a split quaternionic contact structure} if $\paths$ admits a split quaternionic contact structure for which the $(\g, P_{12})$ correspondence space is the given Cartan connection $\eta$ on $\G \to \ph$.

\begin{corollary}[Corollary of Proposition \ref{chsubprop}]\label{chsubcor}
Let $(M, E = V\oplus W)$ be an admissible contact path geometry of dimension at least five with $\ct = 0$. Then $\paths$ admits a split quaternionic contact structure if and only if it has vanishing secondary contact torsion.
\end{corollary}

\begin{proof}
The secondary contact torsion is identified with the homogeneity $(2, -1)$ component of the harmonic curvature of the regular, normal $(\g, P_{12})$ Cartan connection associated to the given contact torsion free contact path geometry. By Theorem \ref{harmoniccorrespondence} the vanishing of the secondary contact torsion implies that the given contact path geometry is the correspondence space of a split quaternionic contact structure on the space of contact paths.
\end{proof}

The homogeneity $(1, 2)$ component of the harmonic curvature of a contact path geometry with vanishing contact torsion and vanishing secondary contact torsion is straightforwardly identified with the harmonic curvature component of the split quaternionic contact structure on the space of contact paths.

\subsubsection{Contact Path Geometries determined by Contact Projective Structures}\label{cprojsection}
Here are made a few remarks about the relation between contact projective structures and general contact path geometries. See \cite{Fox-cproj} for further details on contact projective structures. A contact projective structure, $\en$, on a contact manifold, $(M, H)$, is a null projective structure for the conformal class, $[\theta]$, of contact one-forms. When discussing contact projective structures, let lowercase Latin indices run over $\{1, \dots, 2n-2\}$, and lowercase Greek indices run over $\{0, 1, \dots, 2n-2\}$. If a contact one-form, $\theta$, is fixed, a $\theta$-adapted coframe is a coframe $\theta^{\alpha}$ such that $\theta^{0} = \theta$ and $\theta^{i}(\rb) = 0$, where $\rb$ is the Reeb vector field of $\theta$. If $\theta$ is fixed assume fixed a $\theta$-adapted coframe and dual frame, $E_{\alpha}$, such that $E_{0} = \rb$. Write $\omega = d\theta$. Let $\omega^{ij}\omega_{jk} = -\delta_{k}\,^{i}$, and raise and lower indices with $\omega_{ij}$ according to $\gamma^{p}\omega_{pi} = \gamma_{i}$.

Any representative $\nabla \in \en$ satisfies $\nabla_{(i}\theta_{j)} = 0$. Since the equations defining the geodesics of a connection are independent of the torsion of the connection, there are symmetric representatives $\nabla \in \en$. If $\nabla \in \en$ is symmetric, then $\nabla_{i}\theta_{j} = \tfrac{1}{2}\omega_{ij}$. It is straightforward to check that for symmetric $\nabla \in \en$ the completely trace-free part of $\omega^{kp}\nabla_{p}\omega_{ij}$, which is the section, 
\begin{align*}
-\tfrac{1}{2}\tau_{ij}\,^{k} = \omega^{kp}\left(\nabla_{p}\omega_{ij} + \tfrac{1}{3-2n}\omega^{ab}\left(\omega_{ij}\nabla_{p}\omega_{ab} - \omega_{p[i}\nabla_{j]}\omega_{ab} \right)  \right)\in  \Lambda^{2}(H^{\ast})\tensor H,
\end{align*}
depends on neither the choice of $\theta$ nor the choice of $\nabla$. The tensor $\tau_{ij}\,^{k}$ is called the \textbf{contact torsion} of the contact projective structure. From $\nabla_{[i}\omega_{jk]} = 0$ there follows $\tau_{[ijk]} = 0$. Using $\tau_{[ijk]} = 0$ it is easy to check that $\tau = 0$ if and only if $\tau(X, Y) = 0$ for every isotropic $X \wedge Y \in \Gamma(\Lambda^{2}(H))$. By Theorem C of \cite{Fox-cproj} there is associated to each contact projective structure a curvature normalized regular $(\g, P_{1})$ Cartan connection, $\eta$, on a canonically constructed $P_{1}$ principal bundle $\G \to M$ (the bundle of filtered projective symplectic frames in the tractor bundle). In particular the component $\kappa^{(1)}$ of the curvature of $\eta$ is the contact torsion, and $\eta$ is normal if and only if the contact torsion vanishes. The quotient $\G/P_{12}$ is $\ph$, and so viewing the contact projective structure as a contact path geometry amounts to passing to the correspondence space associated to the inclusion $P_{12} \subset P_{1}$.

\begin{proposition}[\cite{Fox-Schwarzian}]\label{ctorsionequiv}
The contact path geometry induced by a contact projective structure is contact torsion free in the sense of Definition \ref{vancontor} if and only if it has vanishing contact torsion as a contact projective structure. 
\end{proposition}

\begin{proof}
If the contact path geometry underlying a contact projective structure has vanishing contact torsion then by Theorem \ref{prolongationtheorem} it is determined by a regular normal $(\g, P_{12})$ Cartan connection which is, by assumption, the correspondence space of a $(\g, P_{1})$ Cartan connection inducing the given contact projective structure. By Theorem 2.7 of \cite{Cap} this $(\g, P_{1})$ Cartan connection is normal, and by the tabulation of the harmonic curvature components made above, it is regular; it induces the given contact projective structure which by Theorem C of \cite{Fox-cproj} must be contact torsion free.

The contact path geometry underlying a contact projective structure with vanishing contact torsion is that contact path geometry induced on $\ph$ by the correspondence space associated by the inclusion $P_{12} \subset P_{1}$ to the regular normal Cartan connection associated by Theorem C of \cite{Fox-cproj} to the contact projective structure. The correspondence space of a normal Cartan connection is normal, and in the present situation the harmonic part of the curvature function of the correspondence space satisfies $i(x)\kappah = 0$ for all $x \in \p_{1}/\p_{12} = \g_{0, -1}\oplus \g_{0, -2}$. The table of harmonic curvature components above shows that this implies the vanishing of the component of non-positive homogeneity, so that the correspondence space is regular (and normal), and hence induces a contact path geometry with vanishing contact torsion.
\end{proof}

\begin{theorem}\label{vantwotorsiontheorem}
Let $(M, H, [\nabla])$ be an admissible contact projective structure of dimension at least five with vanishing contact torsion. The secondary contact torsion of the contact projective structure viewed as a contact path geometry vanishes if and only if the contact projective structure is flat. In particular, the space of contact paths of an admissible contact projective structure of dimension at least five with $\tau = 0$ admits a codimension $3$ multicontact structure, $C$, so that the canonical projection $\nu:(\ph, \partial(U, W)) \to (\paths, C)$ is a submersion in the category of smooth manifolds equipped with a distribution (i.e.. a split quaternionic contact structure) if and only if the contact projective structure is flat. 
\end{theorem}

\begin{proof}
If $\kappa$ is the curvature function of the correspondence space of the given contact projective structure, Theorem \ref{harmoniccorrespondence} implies that $i(x)\kappah = 0$ for all $x \in \g_{0, -1}$. From the tables given above it is evident that the only possibly non-zero harmonic curvature component is the one of homogeneity $(2, -1)$, which is identified with the harmonic curvature component of the regular normal $(\g, P_{1})$ Cartan connection associated to the underlying contact projective structure. This implies the first claim of the theorem, and the second follows immediately from Corollary \ref{chsubcor}.
\end{proof}

\subsubsection{Totally Geodesic Submanifolds}

An isotropic submanifold, $N \subset M$, is \textbf{totally geodesic} if every contact path tangent to $N$ at any point lies on $N$. Equivalently, $N$ is totally geodesic if $W$ is tangent to $\proj(TN)$.

\begin{example}
  For $3 \leq k \leq n$ the family of $(k-1)$-dimensional isotropic linear submanifolds of $Q_{1}$ parameterized by $Q_{k}$ comprises submanifolds totally geodesic with respect to the flat model contact path geometry.
\end{example}

\begin{proposition}\label{totgeodesicprop}
For a contact path geometry with vanishing contact torsion, if $N \subset M$ is an isotropic submanifold the skew complement in $E^{\perp}$ of $T\proj(TN)$ is well-defined and $N$ is totally geodesic if and only if $T\proj(TN)^{\perp}$ is contained in the restriction to $\proj(TN)$ of $\partial(U, W)$. In particular, a Legendrian submanifold, $N \subset M$, is totally geodesic if and only if $\proj(TN)$ is tangent to $\partial(U, W)$.
\end{proposition}

\begin{proof}
Let $\alpha = r\theta^{-2, -2} + s\theta^{-1, -2}$ with $s \neq 0$ and $\beta = d\alpha$. First it is proved that if $N$ is isotropic then $\proj(TN)$ is $\beta$-null and the skew complement $T\proj(TN)^{\perp}$ in $E^{\perp}$ is well-defined independently of the choice of $\beta$. Because $\theta^{-2, -2}$ is the pullback to $\ph$ of a contact one-form on $M$ and $N$ is isotropic, the restriction to $N$ of $d\theta^{-2, -2}$ vanishes. It follows that the restriction to $N$ of $\beta$ is a non-vanishing multiple of the restriction to $N$ of $d\theta^{-1, -2}$. This shows that the skew complement of $T\proj(TN)$ in $E^{\perp}$ is well-defined independently of the choice of $\beta$. If $A$ and $B$ are tangent to $\proj(TN)$ then so is $[A, B]$, and since $\theta^{-1, -2}$ annihilates $T\proj(TN)$, it follows that $d\theta^{-1, -2}(A, B) = 0$ along $\proj(TN)$, so $\beta(A, B) = 0$ along $\proj(TN)$.

If $N$ is totally geodesic then $W \subset T\proj(TN)$, which by Proposition \ref{twosymp} is equivalent to $T\proj(TN)^{\perp} \subset W^{\perp} = \partial(U, W)$. As $N$ is isotropic, this implies $T\proj(TN) \subset T\proj(TN)^{\perp} \subset \partial(U, W)$. Conversely, if $T\proj(TN)^{\perp} \subset \partial(U, W)$, then Proposition \ref{twosymp} implies $W \subset T\proj(TN)$, so that $N$ is totally geodesic. If $N$ is Legendrian, then $T\proj(TN)$ is a $(2n-3)$-dimensional subbundle of $E^{\perp}$, so $T\proj(TN) = T\proj(TN)^{\perp}$, and hence $N$ is totally geodesic if and only if $\proj(TN)$ is integral for $\partial(U, W)$.
\end{proof}

Proposition \ref{totgeodesicprop} implies that if $N$ is a Legendrian submanifold of an admissible contact path geometry with vanishing contact torsion and secondary contact torsion, then the image in $\paths$ of $\proj(TN)$ is a $(2n-3)$-dimensional submanifold tangent to $C$.

%\bibliographystyle{amsplain}
%\bibliography{cpath}
%Following created using bibtex

\def\cprime{$'$} \def\cprime{$'$} \def\cprime{$'$} \def\cprime{$'$}
  \def\cprime{$'$} \def\cprime{$'$} \def\cprime{$'$} \def\cprime{$'$}
  \def\cprime{$'$} \def\cprime{$'$} \def\cprime{$'$} \def\cprime{$'$}
  \def\cprime{$'$} \def\cprime{$'$} \def\cprime{$'$} \def\cprime{$'$}
  \def\cprime{$'$}
\providecommand{\bysame}{\leavevmode\hbox to3em{\hrulefill}\thinspace}
\providecommand{\MR}{\relax\ifhmode\unskip\space\fi MR }
\providecommand{\MRhref}[2]{%
  \href{http://www.ams.org/mathscinet-getitem?mr=#1}{#2}
}
\providecommand{\href}[2]{#2}

\end{document}